\documentclass[11pt]{amsart}

\usepackage[latin1]{inputenc}
\usepackage[a4paper]{geometry}
\usepackage{amsmath,amsfonts,amssymb,mathrsfs}
\usepackage{graphicx,color}
\def\mk{\medskip}
\def\sk{\smallskip}

\newtheorem{thm}{Theorem}
\newtheorem{prop}[thm]{Proposition}
\newtheorem{cor}[thm]{Corollary}
\newtheorem{lemma}[thm]{Lemma}
\newtheorem{dfn}{Definition}
\newtheorem{rk}{Remark}
\newcommand{\Proof}[1][]{{\par\smallskip\noindent{\bf Proof #1.\enspace }}}
\newcommand{\cqfd}{\hfill\rule{0.35em}{0.35em}}
\newcommand\zu{[0,1)}
\newcommand{\norme}[2][]{\left\|#2\right\|_{#1}} 	
\newcommand\ep{\varepsilon}
\newcommand{\R}{\mathbb{R}}					
\newcommand{\Z}{\mathbb{Z}}					
\newcommand{\Hn}{\mathbb{H}}				
\newcommand{\hn}{\mathfrak{h}}				

\begin{document}

\title[Multifractal analysis on   Heisenberg and Carnot groups]{Multifractal Analysis of functions \\on Heisenberg {and Carnot} Groups}
\author{S.~Seuret}
\author{F.~Vigneron}
\email{stephane.seuret@u-pec.fr}
\email{francois.vigneron@u-pec.fr}
\address{Universit\'e Paris-Est, LAMA (UMR 8050), UPEMLV, UPEC, CNRS, F-94010, Cr\'eteil, France}

\begin{abstract}
{In this article, we investigate the pointwise behaviors of functions on the Heisenberg group. We find wavelet characterizations for the global and local H\"older exponents. Then we prove some a priori upper bounds for the multifractal spectrum of all functions in a given H\"older, Sobolev or Besov space. These upper bounds turn out to be optimal, since in all cases  they are  reached by typical functions in the corresponding functional spaces. We also explain how to adapt  our proof to extend our results to Carnot groups. }
\end{abstract}

\maketitle


\section{Introduction }

 In this article, we draw the first results for a multifractal analysis for functions defined on the Heisenberg group $\Hn$. Multifractal analysis is  now a widespread issue in analysis. Its objective is to provide a description of the variety of local behaviors of  a given  function or a given  measure.  The local behaviors are measured thanks to the pointwise H\"older exponent and one aims at describing the distribution of the iso-H\"older sets \textit{i.e.} the sets of points $x\in \Hn$ with same  pointwise exponent.
 What makes the Heisenberg group interesting for our dimensional considerations is that its Hausdorff dimension is $\dim_H(\mathbb{H}) =4$ while it is defined using only three topological coordinates.  This is due to the special form of the metric in the ``vertical'' direction. This induces surprising properties from the geometric measure theoretic standpoint which are currently being investigated; for instance, Besicovitch's covering theorem and  Marstrand's projection theorem are not true, see \cite {Mat1,Mat2,Besicovitch1,Besicovitch2,Besicovitch3}.  In this paper,  we pursue this investigation by studying the multifractal properties of functions defined on $\Hn$. We find an {\em a priori} upper bound for the Hausdorff dimensions of iso-H\"older sets for all functions in a given  H\"older and Besov space and we prove that these bounds are optimal, since they are reached for {\em generic} functions  (in the sense of Baire's categories) in these function spaces. To do so, we develop methods based on wavelets on $\Hn$ \cite{lemarie89}. In the last Section, we explain how to adapt  our proof to extend our results to   Carnot groups.

 \sk

Let us start by some basic facts on $\Hn$. The first Heisenberg group $\Hn$ consists  \cite[p.~530]{stein} of the set $\R^3$
equipped with a non-commutative group law
$$
(p,q,r)\ast(p',q',r') = (p+p',q+q',r+r'+2(qp'-pq'))
$$
which is also denoted by the absence of multiplicative symbol when the context is clear.
The inverse of $x=(p,q,r)$ is $x^{-1}=(-p,-q,-r)$.
The Haar measure is $dx=dp\wedge dq \wedge dr$ and is also denoted by $\ell$.
It is a homogeneous group \cite[p.~618]{stein} whose dilations are defined by:
$$
\lambda\circ(p,q,r)= (\lambda p,\lambda q,\lambda^2 r).
$$
A left-invariant distance $\delta(x,y)=\norme[\Hn]{x^{-1}  \ast y }$ is given by
the homogeneous pseudo-norm:
\begin{equation}
\label{defmetric}
\norme[\Hn]{x}=\left \{ (p^2+q^2)^2+r^2\right\}^{1/4}.
\end{equation}

\begin{rk}
One can modify the coefficients in formula \eqref{defmetric} so that it becomes a norm on $\Hn$, see \cite{Cygan}. We nonetheless keep this formulation so that the generalization of our results to Carnot groups (see Section \ref{secappendix})  {will require fewer modifications}.
\end{rk}

The Lie Algebra $\hn$ is the vector-space of left invariant
vector fields on $\Hn$. It is nilpotent of step 2 (see \cite[p.~544]{stein}) \textit{i.e.}  $\hn = \mathfrak{n}_1\oplus\mathfrak{n}_2$, where $\mathfrak{n}_1$ is spanned by
\begin{equation}\label{XY}
X=\frac{\partial}{\partial p}+2 q \frac{\partial}{\partial r} \qquad\text{and}\qquad Y=\frac{\partial}{\partial q}-2 p \frac{\partial}{\partial r}
\end{equation}
and $\mathfrak{n}_2$ is spanned by $Z=\frac{\partial}{\partial r}\cdotp$ All commutators
vanish except $[X,Y] = -4 Z$.
Hence $[\hn,\hn]=\mathfrak{n}_2$ and $[\mathfrak{n}_2,\hn]=0$.
The homogeneous structure of $\Hn$ induces dilations of $\hn$:
$$
\lambda \circ\left(\alpha X+\beta Y +\gamma Z\right) = 
\lambda(\alpha X+\beta Y) + \lambda^2 \gamma Z,
$$
which satisfy $\lambda\circ[U,V]=[\lambda\circ U,\lambda\circ V]$.

The positive self-adjoint hypoelliptic Laplace operator on $\Hn$ is (see \cite{folland}):
\begin{equation}
\label{deflaplacien}
\mathcal{L}= -(X^2+Y^2).
\end{equation}
Sobolev spaces of regularity index  $s\geq 0$ can be defined by functional calculus:
\begin{equation}
\label{defsobo}
H^s(\Hn) = \{u\in L^2(\Hn)\,:\, \mathcal{L}^{s/2} u \in L^2(\Hn)\}.
\end{equation}

Throughout the article, $Q=4$ denotes the homogeneous dimension of $\Hn$. In order to   state the results quickly, we postpone  the   classical definitions and notations (horizontal paths, Carnot balls, polynomials, Hausdorff dimension, Besov and H\"older spaces) to \S\ref{sec:def}.

\mk

Let us  define the pointwise H\"older regularity of a function.
\begin{dfn}
Let $f:\Hn\to\R$  be a function belonging to $L^{\infty}_{loc}(\Hn)$. For $s>0$ and $x_0\in\Hn$, $f$ is said to belong to $C^s(x_0)$ if
there exist constants $C>0$,  $\eta>0$ and a polynomial~$P$ with homogeneous degree $\operatorname{deg}_{\Hn}(P)<s$ such that
\begin{equation}\label{Ch}
\forall x\in\Hn, \qquad \|x\|<\eta \quad\Longrightarrow\quad \ |f(x_0 x)-P(x)|\leq C \norme[\Hn]{x}^s.
\end{equation}
One says that $f\in C^s_\text{log}(x_0)$ if, instead of   \eqref{Ch}, the following holds:
\begin{equation}
\label{pertelog} 
|f(x_0 x)-P(x)|\leq C \norme[\Hn]{x}^s \cdot \left|\log\norme[\Hn]{x}\right|
\end{equation}
\end{dfn}

 Observe  that this definition is left-invariant: $f\in C^s(x_0)$ if  and only if
$f_y\in C^s(y^{-1} x_0)$ with $f_{y}:x\mapsto f(y x)$.
 
 \mk
 
The following quantities are crucial in  multifractal analysis.
\begin{dfn}
Let $f:\Hn\to\R$  be a function belonging to $L^{\infty}_{loc}(\Hn)$, and let $x_0\in \Hn$.

The pointwise regularity exponent of $f$ at $x_0$ is
\begin{equation}\label{def_hf}
h_f(x_0)= \sup\{s>0 \,:\, f \in C^s(x_0)\}
\end{equation}
with the convention that $h_f(x_0)=0$ if $f\not\in C^s(x_0)$ for any $s>0$.

\sk
The  multifractal spectrum of $f$ is the mapping $d_f: [0,\infty]\to \{-\infty\}\cup[0,Q]$
\begin{equation*}
d_f(h) = \operatorname{dim}_{H} \left( E_f(h) \right) \quad \text{ where }\quad E_f(h)=\{x\in\Hn \,:\, h_f(x)=h\},
\end{equation*}
where  $\operatorname{dim}_{H}$ stands for the  Hausdorff dimension on $\Hn$. By convention, ${\dim_H\emptyset =-\infty}$. 
\end{dfn}

The multifractal spectrum of  $f$ describes the geometrical distribution of the singula\-rities of $f$ over $\Hn$. The Hausdorff dimension is the right notion to use here, since (at least intuitively, but also for generic functions) the iso-H\"older sets $E_f(h)$ are dense over the support of $f$ and the Minkowski dimension does not distinguish dense sets.

 Wavelets are a key tool in our analysis.
The construction of wavelets on stratified Lie groups has been achieved in \cite{lemarie89}.
A convenient observation is that $\mathcal{Z}=\Z^{3}$
is a sub-group of $\Hn$. For $j\in\mathbb{Z}$ and $k=(k_p,k_q,k_r)\in\mathcal{Z}$, one defines:
$$
x_{j,k}=2^{-j}\circ k = (2^{-j}k_p,2^{-j}k_q,2^{-2j}k_r).
$$
Note that $x_{j,k}^{-1}=x_{j,-k}$.
The dyadic cubes are defined in the following way:
 \begin{gather}\nonumber
 C_{0} = \{(p,q,r) \in \mathbb{H}:  0\leq p,q,r<1\} \ \  \  \mbox{ and } \ \  \
 C_{j,k}  =    x_{j,k}  \ast (2^{-j}\circ  C_{0}).
\end{gather}
\begin{figure}[!h]
\includegraphics[width=280pt,height=140pt]{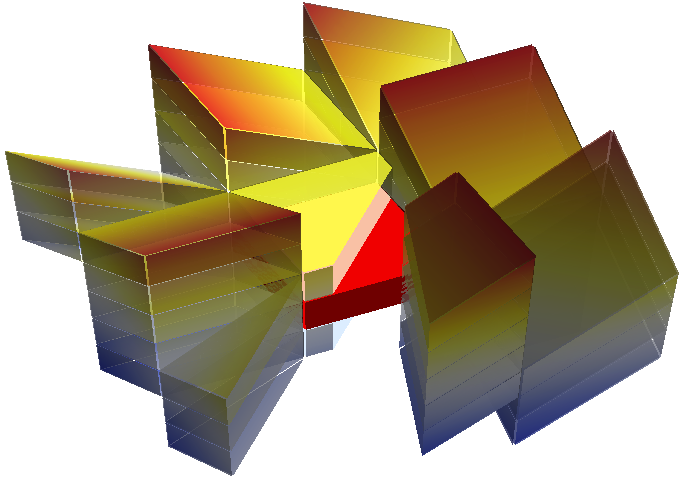}
\caption{\label{CUBES} The cube $C_{0}$ admits $34$ closest neighbors in $\Hn$   contrary to euclidian  cubes of $\R^3$ that admit only $26$ neighbors.}
\end{figure}

The left-multiplication by $x_{j,k}$ maps affine  planes of $\R^3$ on affine  planes, thus the shape of $C_{j,k}$ is a regular parallelogram with vertices on $2^{-j}\circ\mathcal{Z}$. But one shall observe that two different cubes $C_{j,k} $ and $C_{j,k'}$ are not in general euclidian translates of each other.
A neighborhood $\Lambda_{j,k}$ of $C_{j,k}$ is given by 
\begin{equation}\label{DyadicNeighbor}
\Lambda_{j,k}=\bigcup_{k'\in\,\Xi} C_{j,k \ast k'},
\end{equation}
where $\Xi$ is the set of  35 multi-integers $k'=(k_p',k_q',k_r')$ given by   (see Figure~\ref{CUBES}):
$$\begin{array}{ccl}\hline
k_p' & k_q' & k_r' \\\hline
0 & 0 & -1,0,1\\
1 & 0 & -3,-2, -1,0,1\\
1 & 1 & -1,0,1\\
\hline
\end{array}\qquad
\begin{array}{ccl}\hline
k_p' & k_q' & k_r' \\\hline
0 & 1 & -1,0,1,2,3\\
-1 & 1 & 1,2,3\\
-1 & 0 & -1,0,1,2,3
\\\hline
\end{array}\quad
\begin{array}{ccl}\hline
k_p' & k_q' & k_r' \\\hline
-1 & -1 & -1,0,1\\
0 & -1 & -3,-2,-1,0,1\\
1 & -1& -3,-2,-1.\\
\hline
\end{array}
$$ 

\sk
Given $x\in\Hn$ and $j\in\mathbb{Z}$, there exists a unique $k\in\mathcal{Z}$ such that $x\in C_{j,k}$. For this choice of $k$, it is convenient to write
\begin{equation}\label{boxnumeration}
C_j(x)=C_{j,k} \qquad\text{and}\qquad \Lambda_j(x)=\Lambda_{j,k}.
\end{equation}

The diameter of $C_{j,k}$ is $13^{1/4}\times 2^{-j}<2^{1-j}$ (because the diameter of $C_0$ is  $13^{1/4}$).
In particular, if $\delta(x,y)<2^{-j}$, then $x,y$ belong simultaneously to at least one $\Lambda_{j,k}$.

\mk

 Let us recall now the construction of wavelets on $\Hn$ by Lemari\'e   \cite{lemarie89}. 
For any integer $M>Q/2$, there exist $2^Q-1=15$ functions $(\vartheta_\ep)_{ 1\leq \varepsilon \leq  15} $ in $ H^{4M}(\Hn)$ such that:
  
\begin{itemize}
\item  
There exist $C_0,r_0>0$ such that for any multi-index $\alpha$  of
length $|\alpha|<4M-Q$:
\begin{equation}\label{DECAYPsi}
\forall x\in\Hn, \qquad
|\nabla_{\Hn}^\alpha \vartheta_{\varepsilon}(x)| \leq C_0 \exp\left(- {\norme[\Hn]{x}}/{r_0}\right).
\end{equation}

\sk

\item
Each function $\Psi_{\varepsilon}= \mathcal{L}^M \vartheta_{\varepsilon}$ has $2M$ vanishing moments, i.e. for every   polynomial function $P$ of homogeneous degree $\operatorname{deg}_{\Hn}P<2M$,
then 
\begin{equation}
\label{vanishingmoments}
\int_{\Hn} \Psi_\varepsilon (x)P(x)dx=0.
\end{equation}
Moreover,  $|\Psi^\varepsilon (x)|\leq C_0 \exp\left(-\norme[\Hn]{x}/r_0\right)$ and
$\Psi^\varepsilon  \in H^\sigma(\Hn)$ for $\sigma<2M-Q$.

\sk
\item The family of functions $(2^{jQ/2}\Psi^\varepsilon_{j,k})_{j\in Z,k\in \mathcal{Z},  1\leq \varepsilon \leq  15}$, where    
 $$
\Psi^\varepsilon_{j,k}(x) = \Psi_\varepsilon\left(2^j\circ (x_{j,k}^{-1} \ast x)\right),
$$
forms a  Hilbert basis of $L^2(\Hn)$, \textit{i.e.}:
\begin{equation}\label{WAVELETS}
f  \, \stackrel{ \text{$L^2$} }{=} \,  \sum_{\varepsilon,j,k} d^\varepsilon_{j,k}(f)\Psi^\varepsilon_{j,k}
\quad\text{with}\quad
 d^\varepsilon_{j,k}(f) = 2^{jQ}\int_{\Hn} f(x)\Psi^\varepsilon_{j,k}(x)dx.
\end{equation}

\end{itemize} 
 
The real numbers $ d^\varepsilon_{j,k}(f)$ are called the {\em wavelet coefficients} of $f$. Note that we use an $L^\infty$ normalization for  the wavelet  in \eqref{WAVELETS} and that our choice implies   that  the family $(2^{-jQ/2}d_{j,k}^\varepsilon(f))$ belongs to $\ell^2 $ and thus   tends to 0 when $j\to\pm\infty$ and $\norme[\Hn]{k}\to\infty$.

\begin{rk}
Instead of Lemari\'e's wavelet, one may use the wavelets built by F\"uhr and Mayeli in \cite{FuhrMayeli} which can be constructed on any Carnot group. We will use these wavelets when explaining the generalization of our results to more general Carnot groups.
\end{rk}

\mk
We can now state our main theorems. 
For non-integer regularity, H\"older classes   can be totally described with wavelets coefficients, as in the Euclidian case. 
\begin{thm}
\label{thglobal}
For $s\in\R_+\backslash\mathbb{N}$  and $[s]<2M$,
a function $f$  belongs to $ C^s (\Hn)$ if and only if there exists a constant $C>0$ such that
\begin{equation}\label{GLOBALHOLDER}
\forall (\varepsilon,j,k)\in\{1,\ldots,2^Q-1\}\times\Z\times\mathcal{Z},\qquad |d^\varepsilon_{j,k}(f)| \leq  C 2^{-js}.
\end{equation}
\end{thm}
Theorem \ref{thglobal} is essentially proved in \cite{FuhrMayeli}, Theorems 5.4 and 6.1. We give another proof here, using    Lemari\'e's wavelets. The  existence of the decomposition \eqref{WAVELETS} allows us to obtain a straightforward equivalence between $C^s(\Hn)$ and \eqref{GLOBALHOLDER}.

\mk
Up to a logarithmic factor, the pointwise regularity class $C^s (x_0)$ can also
be described with wavelets coefficients.
\begin{thm}\label{thlocal}
Given $f\in L^2(\Hn)$ and $x_0\in \Hn$, the following properties hold.
\begin{itemize}
\item
If $f\in C^s(x_0)$, then there is $R>0$ such that for any indices $\varepsilon,j,k$:
\begin{equation}\label{PointRegWav}
\delta(x_{j,k},x_0)<R \quad\Longrightarrow\quad
|d^\varepsilon_{j,k}(f)| \leq  C 2^{-js} \left(1+2^j \delta(x_{j,k},x_0)\right)^s.
\end{equation}
\item
Conversely, if $f$ satisfies~\eqref{PointRegWav}  and belongs to $C^{\sigma}(\Hn)$ for an arbitrary small $\sigma>0$, then~$f$ belongs  to $C^s_{\log}(x_0)$.
\end{itemize}
\end{thm}
 
 \begin{rk} The important information contained in \eqref{PointRegWav} does not just lie in the coefficients closest (at each dyadic scale) to $x_0$:
\begin{equation*}
|d_{j,k}^\varepsilon (f)| \leq \begin{cases}
C2^{-js} &\text{if }\delta(x_{j,k},x_0)\lesssim 2^{-j}\\
C\delta(x_{j,k},x_0)^s &\text{if } 2^{-j}\lesssim\delta(x_{j,k},x_0)< R.
\end{cases}
\end{equation*}
\end{rk} 
\begin{rk} In the Euclidian case, wavelet leaders  \cite{jaffleaders} are more stable  numerically.
The wavelet leaders of a function $f\in L^2(\Hn)$ is the sequence
\begin{equation*} 
D_{j}(f,x) = \sup\left\{|d^\ep_{j',k'}(f)| : j'\geq j \mbox{ and } C_{j',k'}\subset \Lambda_j(x)\right\}
\end{equation*}
with $\Lambda_j(x)$ defined by~\eqref{boxnumeration}.
One checks easily that another statement  equivalent to \eqref{PointRegWav}~is:
\begin{equation*}
f\in C^s(x_0) \quad\Longrightarrow\quad 
\forall \ j\geq 0, \ \ D_j(f,x_0)\leq C 2^{-js}.
\end{equation*} 
\end{rk}

It is also obvious from the last two theorems that $f\in C^s(\Hn)$ implies  $h_f(x)\geq s$ for every $x\in \Hn$. The optimality of this result is asserted by the following theorem. Recall that a property~$\mathcal{P}$ is \textit{generic}   in a complete metric space $E$ when it holds on 
a \textit{residual} set \textit{i.e.} a set  with a complement of first Baire category. A set is of first Baire category if it is the 
union of countably many nowhere dense sets. As it is often the case,  it is enough to build a 
residual set  which is  a countable  intersection of dense open sets in $E$. 

\begin{thm}
\label{thgenericCalpha}
There exists a dense open set (hence a generic set) $\mathcal{R}$ of functions    in  $ C^s(\Hn)$ such that for every  $f\in \mathcal{R}$
and  every $x\in \Hn$, $h_f(x)=s$.   
\end{thm}
\noindent
In particular, generic functions   in  $ C^s(\Hn)$  are monofractal \textit{i.e.} $E_f(h)=\emptyset $ if $h\neq s$.
 
 \mk 
 
One can also obtain  {\em a priori} upper bounds for the multifractal spectrum of functions belonging to Besov and Sobolev spaces on $\Hn$ (see Section \ref{secdefsobo} for precise definitions).
These function spaces play a fundamental role in harmonic and functional analysis.
There are two types of results one can naturally look for: general local regularity results that are valid for all functions in a given Besov space, and results that are only true for ``almost every'' function
in this function space.

\begin{thm}
\label{majspectrumBesov}
For $s> {Q}/{p}$, every $f\in B^s_{p,q}(\Hn)$ satisfies:
\begin{equation}\label{BesovSpectrum}
\mbox{for all } h\geq s-Q/p, \ \ \ d_f(h)\leq \min\left(Q \,,\, p\left(h-s+ {Q}/{p}\right)\right)
\end{equation}
and $d_f(h)=-\infty$ if $h<s-Q/p$.
 
\end{thm}

This theorem has many remarkable consequences. For instance, it illustrates the optimality of the Sobolev inclusion  $B^{s}_{p,q}(\Hn) \hookrightarrow C^{s-Q/p}(\Hn)$ : the sets of points with the least possible pointwise H\"older exponent $s- Q/p$ has Hausdorff dimension at most~0. Similarly, as a consequence of the proof, the set of points whose pointwise H\"older exponent  is at least~$s$ and has  full Haar  measure in $\Hn$.
The main difference with $C^s(\Hn)$ is that functions in $B^s_{p,q}(\Hn)$ may really be multifractal, meaning that many iso-H\"older sets $E_f(h)$ are non-empty with a non-trivial Hausdorff dimension.  This is the case for generic functions in $B^s_{p,q}(\Hn)$.

\begin{figure}[]
\includegraphics[width=160pt]{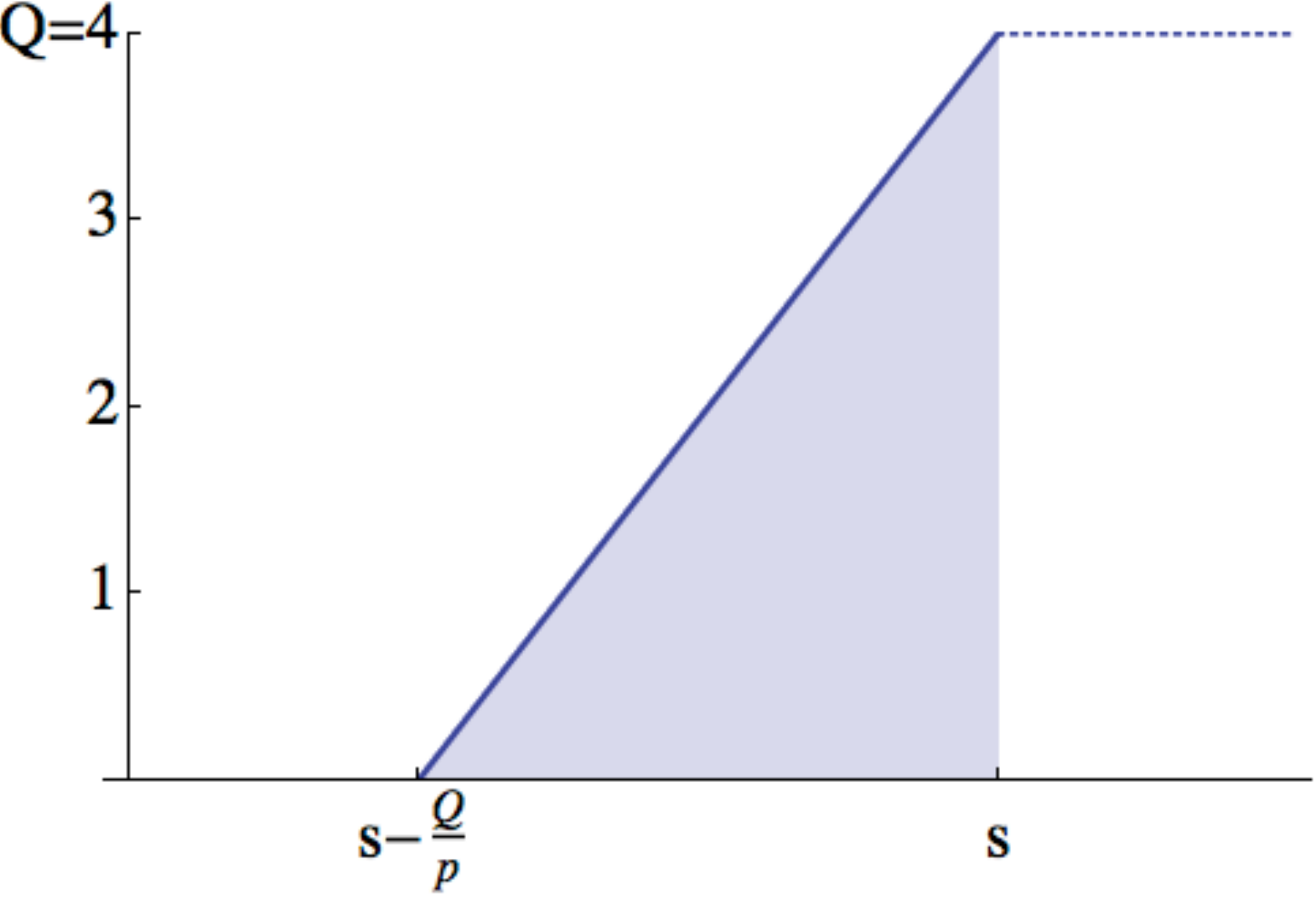}
\caption{Upper bound for  the multifractal spectrum of functions in
 $B^s_{p,q}(\Hn)$.}
\end{figure}

\begin{thm}\label{th:besovps}
For $s> Q/p$, there is a residual set $\widetilde{ \mathcal{R}} \subset B^s_{p,q}(\Hn)$ such that for all $f \in \widetilde{ \mathcal{R}}$,  
\begin{equation}\label{BesovSpectrumequality}
\forall h\in [s- Q/p  ,s], \quad d_f(h)=   p\left(h-s+ {Q}/{p}\right),
\end{equation}
and  $E_f(h) = \emptyset $ for all other exponents.
\end{thm} 
\noindent
In particular, for generic functions $f \in B^s_{p,q}(\Hn)$, Haar-almost every point has a pointwise H\"older exponent equal to $s$.

\mk

This paper is organized as follows. Section \ref{sec:def} contains   the  definitions and previous results that we use in the sequel. In  Sections  \ref{secglobal} and  \ref{secpointwise} respectively, we deal with global and pointwise H\"older regularity (Theorems \ref{thglobal} and \ref{thlocal}). In particular, the monofractality of generic functions  (Theorem  \ref{thgenericCalpha}) in $C^s(\Hn)$ is proved in Section \ref{genericcalpha}. The multifractal properties of functions in a Besov space are then investigated in Sections \ref{secbesov1} and  \ref{section:generic_besov}. Finally, we explain how to extend our results to general stratified nilpotent groups (\textit{i.e.}~Carnot groups) in Section \ref{secappendix}.

\mk

Let us finish this introduction with a question. It would be very interesting to be able to represent the functions on $\Hn$, or at least the traces of such functions on affine subspaces of $\R^3$. Indeed, the natural anisotropy induced by the metric on $\Hn$ should create an anisotropic picture, and as of today, creating natural and simple models for anisotropic textures is a great challenge in image processing. Of course the starting point would be to understand how to draw a wavelet (be it a Lemari\'e wavelet, or another one!)  on~$\Hn$. We   believe that this is a very promising research direction.


\section{Definitions and recalls}
\label{sec:def}


\subsection{Balls on $\Hn$}

As the shape of balls is rather counter-intuitive on $\Hn$,
a few geometric statements will be useful in the following. The volume of the gauge balls
$$
\mathcal{B}(x, r)=\{y\in\Hn\,: \,\delta(x,y)<r\}
$$
is denoted by $\ell(\mathcal{B}(x,r))$ and is equal to $ \frac{\pi^2}{2}r ^Q$ with $Q=4$. In particular, the Haar measure $\ell$ has
the doubling property:  $l(B(x,2r)) \leq C l(B(x,r))$ for some universal constant $C$.

\mk
  For any $x,x'\in\Hn$ and $r,r'>0$, one has $
x'\ast(r'\circ\mathcal{B}(x,r))=\mathcal{B}(x'\ast(r'\circ x), r r')$ and  $\mathcal{B}(x,r) = x\ast (r \circ \mathcal{B}(0,1))$.

\mk
The triangular inequality holds   in general with a constant depending on the metric.

 \begin{prop}[Folland, Stein, prop.1.6]
There exists a constant $\gamma_1>0$ such that
\begin{equation}\label{TRIANGULAR}
\forall x,y\in\Hn, \quad \norme[\Hn]{xy} \leq \gamma_1 \left(  \norme[\Hn]{x}+ \norme[\Hn]{y} \right).
\end{equation}
\end{prop}
\noindent In particular, the diameter of a gauge ball $\mathcal{B}(x,r)$ does not exceed $2\gamma_1 r$. One will use this property later in the following form:
\begin{cor}\label{CorTRIANGULAR}
There exists $C>0$ such that for any $\eta>0$ and $x\in \mathcal{B}(0,\eta/C)$, one has $$\mathcal{B}(0,\eta/C)\subset \mathcal{B}(x,\eta).$$
\end{cor}

\subsection{Polynomial functions}

A polynomial function $P$ on $\Hn$ is a polynomial function of the coordinates $(p,q,r)$ ; its homogeneous degree is defined by
$$\deg_\Hn P = \deg P(t,t,t^2)$$
where the right-hand side is computed in $\R[t]$.
Given $x=(p,q,r)$ and $\alpha\in\{1,2,3\}^m$, one defines
$x^{\alpha} = x(\alpha_1)\ldots x(\alpha_m)$
with $x(1)=p$, $x(2)=q$ and $x(3)=r$.
A polynomial function $\Hn$ of homogeneous degree at most $N$ is thus a function of the form:
$$
P(x)=\sum_{|\alpha|\leq N} c_\alpha x^\alpha
$$
with $c_\alpha\in\R$ and $|\alpha|=\sum \omega(\alpha_i)$ with $\omega(1)=\omega(2)=1$ and $\omega(3)=2$.

\subsection{The operator $\nabla_\Hn$}

Let us denote by $\nabla_\Hn=(X,Y)$ the basis  (equation \eqref{XY}) of horizontal derivatives.
Given a multi-index $\alpha\in\{1,2,3\}^m$ one will denote by
\begin{equation}
\label{nabla_Hn}
\nabla_\Hn^\alpha f = V_{\alpha_1}\ldots V_{\alpha_m}f
\end{equation}
where $V_1=X$, $V_2=Y$ and $V_3=Z$. As $Z=-\frac{1}{4}[X,Y]$, one may reduce $\nabla_\Hn^\alpha f$ to a linear combination of  terms that contain exactly $|\alpha|$ powers of $X$ and $Y$. One says that $\nabla_\Hn^\alpha f$ is a \textit{horizontal derivative} of $f$ of order $|\alpha|$.

\subsection{Horizontal paths and Taylor formula}

In this short Section, we elaborate on Taylor expansions, see \cite{follandstein,Varadarajan} for details. We emphasize the notions needed later in this paper.

In this short Section, we elaborate on Taylor expansions 
Two points $x,y\in \Hn$ can always be joined by a sub-unitary horizontal path, \textsl{i.e.} a piecewise Lipschitz    arc $\gamma:[0,L]\to\Hn$  such that
for almost every $t$, the tangent vector can be decomposed as
$$\gamma'(t) = \alpha(t) X(\gamma(t))+\beta(t) Y(\gamma(t))$$
with $\alpha^2(t)+\beta^2(t)\leq 1$. The so-called Carnot-length $d_C(x,y)=\inf_\gamma \left(\int_0^L \alpha^2(t)+\beta^2(t) dt\right)^{1/2}$ is uniformly equivalent to $\delta(x,y)$.
Integrating along such an arc provides the first order Taylor formula:
$$f(y)=f(x)+\int_0^{L} \nabla_{\Hn}f(\gamma(t)) \gamma'(t) dt.$$
In turn, this identity provides a Lipschitz estimate.
\begin{prop}[\cite{follandstein},  Theorem 1.41]\label{lipschitz}
There exists $C>0$ and $\gamma_2>0$ such that for all $f\in C^1(\Hn)$ and $x,y\in\Hn$,
$$|f(y)-f(x)|\leq C \delta(x,y) \sup_{\norme{z}\leq \gamma_2\delta(x,y)} |\nabla_\Hn f(xz)|.$$
\end{prop}

The left-invariant Taylor expansion of a function is given by the next definition.
\begin{dfn}[\cite{follandstein}]
The right Taylor polynomial of homogeneous degree $k$ of a smooth function~$f$ at $x_0\in\Hn$ is the unique
polynomial $P_{x_0}$ of homogeneous degree $\leq k$ such that
$$\forall\alpha\in\bigcup_{m\in\mathbb{N}}\{1,2,3\}^m, \qquad |\alpha|\leq k \quad\Longrightarrow\quad \nabla_\Hn^\alpha f(x_0) = \nabla_\Hn^\alpha P_{x_0}(0).$$
\end{dfn}

To proceed with the subsequent calculations we will need to write down the Taylor expansion \textsl{explicitly}.
It must be done carefully for various reasons. The most obvious one is that  $XYf\neq YXf$ but~$pq=qp$.
The second problem induced by the   anisotropy of the Heisenberg structure is that
the traditional match between the index of the derivative and the index of the polynomial  breaks down.
For example, at the 2$^{\text{nd}}$ order near the origin, $f(p,q,0)$ will be computed using only the first 2 powers of $p$ and $q$ but, contrary to the euclidian setting,
it will involve vertical derivatives at the origin, through $Zf(0)$. 

With this in mind, a good way to write the Taylor polynomial  of order $N$ down is:
\begin{equation}\label{TAYLOR}
P_{x_0}(y) =  \sum_{k=0,...,N}  \ \sum_{|\alpha|=   k}  \  y^\alpha \ \Big( \sum_{|\beta| =k } c_{\alpha,\beta} \nabla_\Hn^\beta f(x_0) \Big) = \sum_{|\alpha|=  |\beta|\leq N }  c_{\alpha,\beta} \nabla_\Hn^\beta f(x_0) y^\alpha.
\end{equation}
Beyond order 2, even though the polynomial $P_{x_0}$ remains unique, the coefficients $c_{\alpha,\beta}$  in~\eqref{TAYLOR} are not and thus a choice has to be done once for all
before starting a computation. For example, one possible writing of  the polynomial of order 3 at the origin is:
\begin{align*}
P_0(p,q,r) &= f(0)+p Xf(0) + q Yf(0)\\
&+\frac{1}{2}\left(p^2 X^2f(0)+ 2pqXYf(0)+q^2 Y^2f(0)\right)+({\color{blue}2pq}+r)\cdot{\color{blue}Z}f(0)\\
&+\frac{1}{3!}\left(p^3X^3f(0)+3p^2qX^2Yf(0)+3pq^2XY^2f(0)+q^3Y^3f(0)\right)\\
&+({\color{blue}2pq}+r)\cdot\left(pX{\color{blue}Z}f(0)+qY{\color{blue}Z}f(0)\right).
\end{align*}
 The actual choice between the possible  expressions  is irrelevant. Since monomials are commutative and $Z=-\frac{1}{4}[X,Y]$ one can assume from now on that the formula is reduced to  indices   $\beta\in\bigcup_{m\in\mathbb{N}}\{1,2\}^m$.   Further results and explicit Taylor formulas on homogenous groups can be found in  \cite{TaylorBonfiglioli}.

\begin{figure}
\begin{center}
\includegraphics[width=62mm]{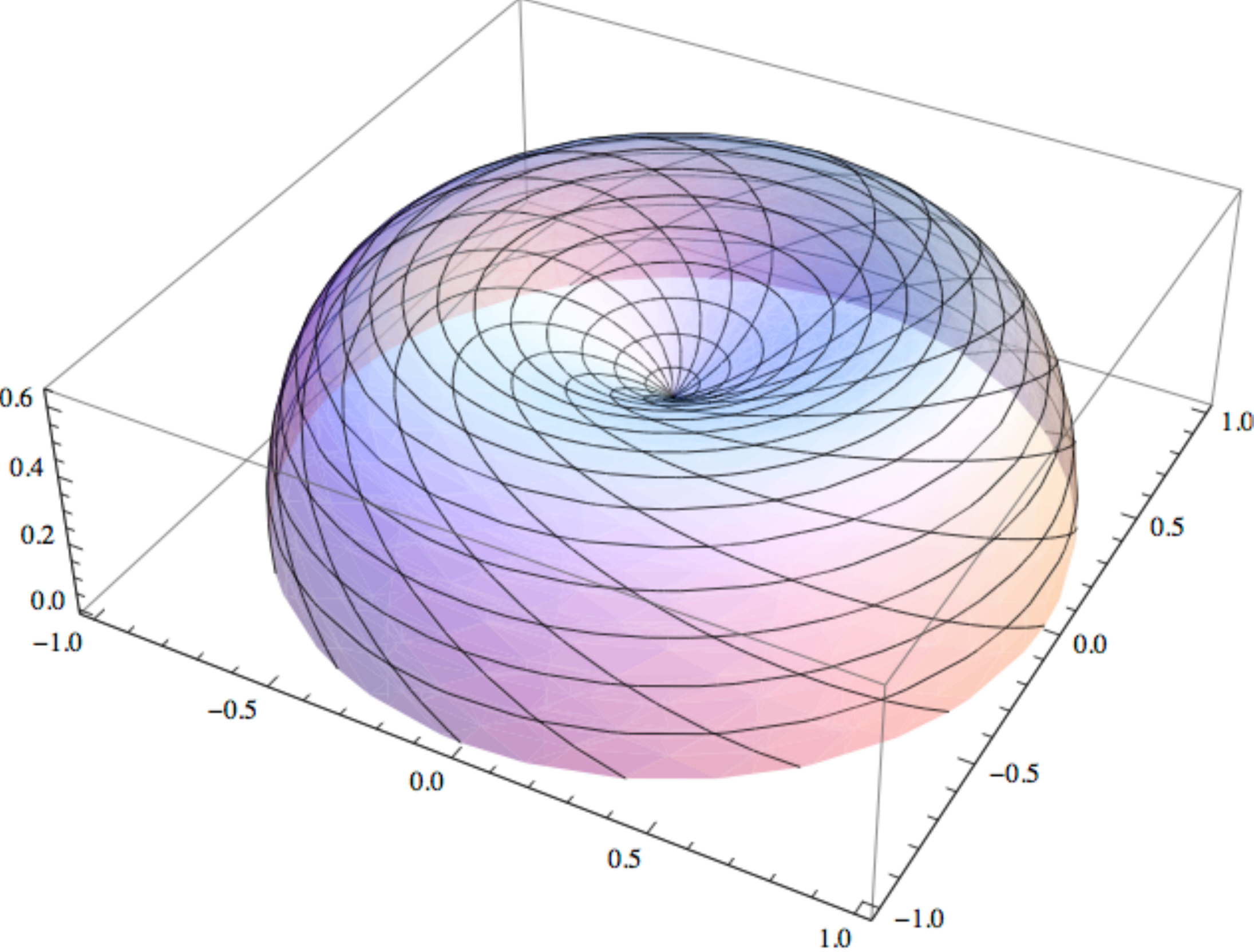}\qquad
\includegraphics[width=40mm]{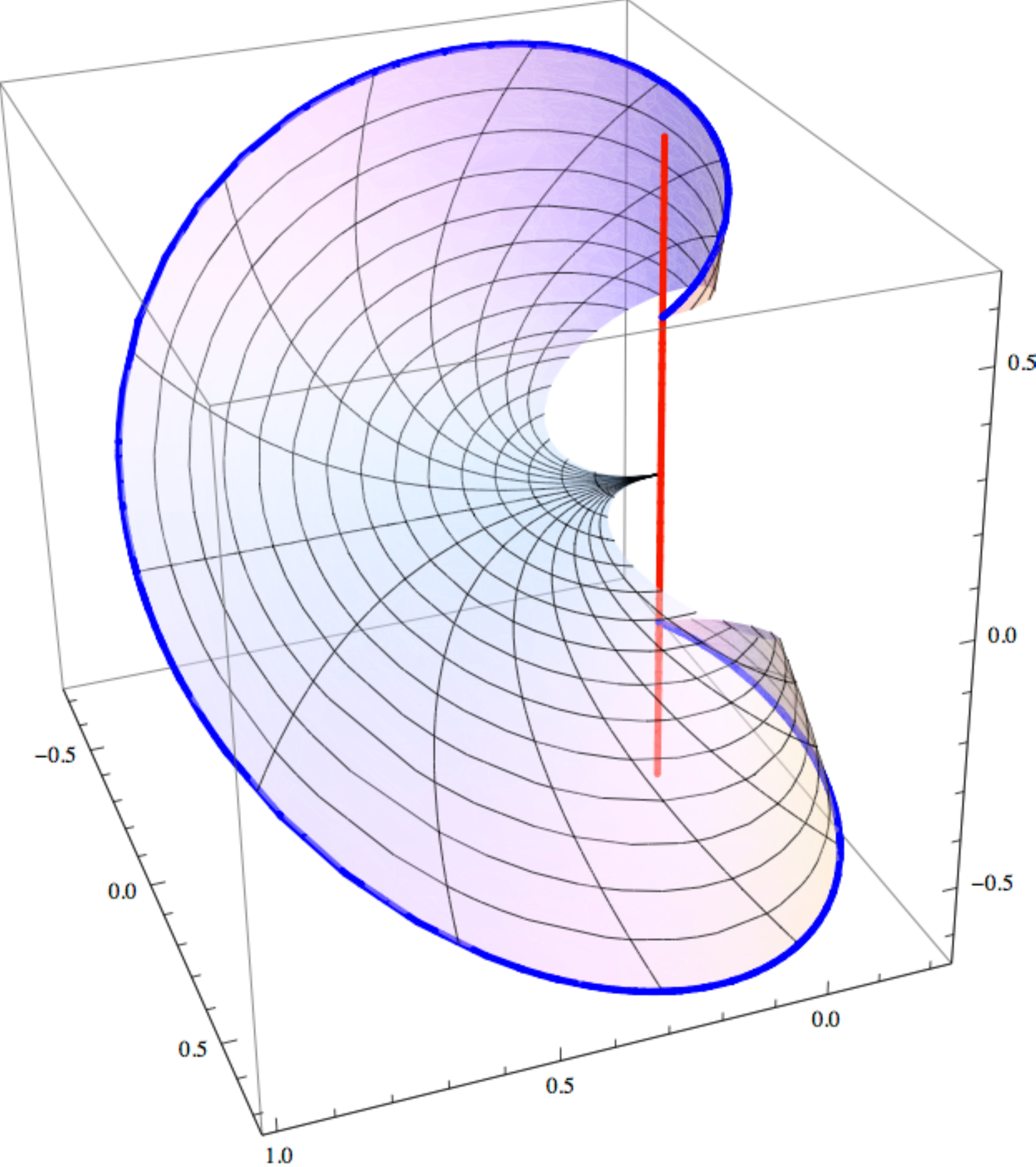}
\caption{Left : Upper half of the unit Carnot ball of $\Hn$. 
 Right: one meridian arc of the unit ball and the
unitary ``horizontal'' geodesics joining the origin to each point of this meridian arc. }
\end{center}
\end{figure}

As expected,  the right Taylor polynomial approximates $f(x_0y)$ for $y$ small enough.

\begin{thm}[Folland, Stein, corr. 1.44]\label{taylor}
If $f\in C^{k+1}(\Hn)$, then the following estimate holds for some universal constant $C_k$:
$$
|f(x_0y)-P_{x_0}(y)| \leq C_k \norme[\Hn]{y}^{k+1}  \sup_{|\alpha|=k+1} \big(\sup_{\norme{z}\leq \gamma_2^{k+1}}|\nabla_\Hn^\alpha f(x_0z)|\big)
$$
\end{thm}

\subsection{Hausdorff dimension on $\Hn$}

The diameter of a set $E\subset \mathbb{H}$ will be denoted by~
$$
|E|=\sup\{\delta(x,y): x,y\in E\}.
$$ 
Let us recall the definition of  the Hausdorff measures and dimension. Let $s>0$ and $\eta>0$ be two positive real numbers. For any set $E\subset \mathbb{H}$, one defines
\begin{equation}
\label{defHmeasure}
\mathcal{H}^s_\eta (E) = \inf _{\mathcal{R}}  \ \sum_{B\in\mathcal{R}} |B|^s \ \ \in [0,+\infty],
\end{equation}
where the infimum is taken over all possible coverings $\mathcal {R}$ of $E$ by gauge balls of radii less than~$\eta$. Recall that a covering of $E$ is a family~$\mathcal{R}=\{B_i\}_{i\in I}$ of balls  satisfying $$E \subset \bigcup_{i\in I} B_i  .$$
The mapping $\eta\mapsto \mathcal{H}^s_\eta (E)$ is decreasing with $\eta$,  hence one can define
$$\mathcal{H}^s  (E)= \lim_{\eta\to 0^+} \ \mathcal{H}^s_\eta (E)  \ \ \in [0,+\infty].$$
From this definition, it is standard to see that $s\mapsto \mathcal{H}^s  (E)$ is a decreasing function that jumps from infinity to zero at a unique real number called the Hausdorff dimension of~$E$:
$$
\dim_H E = \inf\{s:\mathcal{H}^s  (E)=0\} = \sup\{s : \mathcal{H}^s  (E)=+\infty\}.
$$

 \subsection{ H\"older, Sobolev and Besov  regularity}
\label{secdefsobo}

 \begin{dfn}
For $s=k+\sigma$ with $k\in\mathbb{N}$ and $\sigma\in]0,1[$,   $C^s(\Hn)$ is the set of   functions such that  for any multi-index of length $|\alpha|\leq k$, the function $\nabla_{\Hn}^\alpha f$ is continuous and:
\begin{equation*}
\sup_{|\alpha|=k}\frac{|\nabla_{\Hn}^\alpha f(x)-\nabla_{\Hn}^\alpha f(y)|}{\delta(x,y)^\sigma}<\infty. 
\end{equation*}
\end{dfn}

 Hence, 
  H\"older classes are  defined as in the Euclidian case (see \cite{CCX,CX}).  Two equivalent Banach norms on $C^s(\Hn)$, denoted by~$\norme[C^s(\Hn)]{f}$, are given by
$$\sup_{|\alpha|=[s]}\frac{|\nabla_{\Hn}^\alpha f(x)-\nabla_{\Hn}^\alpha f(y)|}{\delta(x,y)^{s-[s]}} 
\qquad\text{and}\qquad
\sup_{\varepsilon,j,k} \left(2^{js} |d_{j,k}^\varepsilon(f)|\right).$$

Sobolev spaces,  that we introduced before in \eqref{defsobo}, can also be described using the horizontal derivatives \eqref{nabla_Hn} (see  \cite{MV} and  references therein, and also \cite{folland} for more details on these functional spaces).
\begin{prop} 
For $k\in\mathbb{N}$, a function~$f$ belongs to~$H^k(\Hn)$ if and only if
for any multi-index $\alpha$ of length~$|\alpha|\leq k$,
$
\nabla_{\Hn}^\alpha f \in L^2(\Hn).
$

\sk

For $s=k+\sigma$ with $k\in\mathbb{N}$ and $\sigma\in]0,1[$, one has
$f\in H^s(\Hn)$ if and only if $f\in H^k(\Hn)$ and for $|\alpha|=k$ (with $Q=4$, the homogeneous dimension of $\Hn$):
\begin{equation}\label{EQ:MV}
\iint_{\Hn\times\Hn} \frac{| \nabla_{\Hn}^\alpha f (x)-\nabla_{\Hn}^\alpha f(y)|^2}{\delta(x,y)^{Q+2\sigma}} \:dxdy<\infty.
\end{equation}
\end{prop}

One has the continuous inclusion  $H^s(\Hn)\subset C^{s-Q/2}(\Hn)$, which holds if $s>Q/2$
and $s-Q/2\not\in\mathbb{N}$. Namely, for any multi-index $\alpha$ such that
$s=Q/2+|\alpha|+\sigma$ and $\sigma\in]0,1[$:
$$
|\nabla_{\Hn}^\alpha f(x)-\nabla_{\Hn}^\alpha f(y)|
\leq C_s \norme[H^s(\Hn)]{f} \delta(x,y)^{s-Q/2}.
$$

As in the euclidian case, this inclusion fails
if $s-Q/2\in\mathbb{N}$, because one can then find
$f\in H^s(\Hn)$ and $|\alpha|=s-Q/2$ 
such that $\nabla_{\Hn}^\alpha f$ is not a continuous function (note that $Q=4$ is even).
However, the corresponding inclusion in $\text{BMO}$ holds (see \cite{CX}).

\mk
On $\R^n$, Besov spaces can be defined in various ways and the equivalence
between all those definitions is part of the folklore. For nilpotent Lie groups,
the situation is less straightforward and all the equivalences should be checked
carefully. 
For multifractal analysis, the most convenient definition of Besov spaces is:
\begin{dfn}
The Besov space $B^s_{p,q}(\Hn)$ of \cite{saka} consists of functions $f$ on $\Hn$ such that:
\begin{equation}
\label{defbesov}
a_j = \Big\| 2^{j\left(s- Q/p\right)}d^\varepsilon_{j,k}(f) \Big\|_{\ell^p(k)}\in \ell^q(j).
\end{equation}

\end{dfn}
Other definition of Besov spaces involve    Littlewood-Paley theory   (see \cite{stein_LP} for a historical background) and continuous wavelet decomposition (see \cite{FuhrMayeli} for instance), but all definitions coincide. The interested reader can have a look at Section 5 of \cite{Ajoutreferee}.  
Depending on the applications, other definitions have proved useful.
Trace theory  \cite{DGN,vigneron}  is better
understood with the geometric norm \eqref{EQ:MV}, real interpolation
norms with operator theory  \cite{BP,saka}, while 
complex interpolation norms relate to microlocal analysis and Weyl
calculus on~$\Hn$   \cite{BC,CCX,MV}.


\section{Global H\"older regularity with wavelets coefficients: Theorem \ref{thglobal}}
\label{secglobal}

We first prove in Section \ref{sec31} that the wavelet coefficients of every function $f\in C^s(\Hn)$ enjoy the decay property \eqref{GLOBALHOLDER}. The converse property is shown in Section \ref{sec32}.

 \subsection{Upper bound for the wavelets coefficients}
 \label{sec31}
  Assume that $f\in C^s(\Hn)$.  
Let $s= {[s]}+\sigma$ with $ [s]\in\mathbb{N}$ and $0<\sigma<1$.
The change of variables $y=2^j\circ (x_{j,k}^{-1} x)$ reads:
$$
d_{j,k}^\varepsilon=\int_{\Hn}f(x_{j,k} (2^{-j}\circ y)) \Psi_\varepsilon(y) dy.
$$

When $[s]=0$, one infers from the vanishing moment of $\Psi_\ep$ (i.e. $\int_{\Hn} \Psi_\varepsilon(x)dx =0$) that:
\begin{eqnarray*}\nonumber
|d_{j,k}^\varepsilon|  \!\!  \!\!  & =& \!\!   \!\!   \left|\int_{\Hn} (f(x_{j,k} (2^{-j}\circ y))-f(x_{j,k})) \Psi_\varepsilon(y) dy\right|   \!   \leq \!    \int_{\Hn}\! |f(x_{j,k} (2^{-j}\circ y))-f(x_{j,k})| |\Psi_\varepsilon(y)| dy \\
 \!\!  &\leq& \!\! 
  2^{-js}\norme[C^s(\Hn)]{f} \int_{\Hn} \norme[\Hn]{y}^s |\Psi_\varepsilon(y)| dy = C 2^{-js}.\nonumber
\end{eqnarray*}

\sk

When $[s]\geq1$,  one uses $\Psi^\varepsilon = \mathcal{L}^M\vartheta_\varepsilon$ and proceeds with ${[s]}$ integrations by part against the function $
g_{j,k}(y)=f(x_{j,k}  (2^{-j}\circ y)).
$

 Observe that   the homogeneity of the horizontal derivatives     yields   that for every $\alpha$, 

\begin{equation}\label{DerivativesWithScaling}
\nabla_\Hn^\alpha \,  g_{j,k}(y)= 2^{-j|\alpha|} \times [\nabla_\Hn^\alpha f]((x_{j,k}  \ast (2^{-j}\circ y)).
\end{equation}

When ${[s]}=2m$ is even ($m\geq 1$), one writes  $
d_{j,k}^\varepsilon = \int_{\Hn} (\mathcal{L}^m g_{j,k})(x) \widetilde\vartheta_\varepsilon (x) dx$, 
with $\widetilde\vartheta_\varepsilon = \mathcal{L}^{M-m}\vartheta_\varepsilon$.  The term  $ \mathcal{L}^{m}g_{j,k}=(-X^2-Y^2)^{m}(g_{j,k})$ can be developed, and one gets
$$ \mathcal{L}^{m}g_{j,k}= \sum_{|\alpha|={[s]}}  l_\alpha  \nabla_\Hn^\alpha(g_{j,k}),$$ 
 for some coefficients $l_\alpha$ independent of the problem. Recalling \eqref{DECAYPsi},  $\widetilde\vartheta_\varepsilon$ is a  well-localized function. One can use the vanishing moments of $\vartheta_\varepsilon$ (and thus  of $\widetilde\vartheta_\varepsilon$)  to get
\begin{eqnarray*}
d_{j,k}^\varepsilon   & = &  \sum_{|\alpha|={[s]}}  l_\alpha 
\int_{\Hn}\left(\nabla_\Hn^\alpha g_{j,k}(y)-\nabla_\Hn^\alpha g_{j,k}(0)\right)
\widetilde\vartheta_{\varepsilon }(y)dy, \\
& = & 2^{-j[s]}  \sum_{|\alpha|={[s]}}  l_\alpha 
\int_{\Hn}    \left(\nabla_\Hn^\alpha f((x_{j,k} \ast (2^{-j}\circ y))-\nabla_\Hn^\alpha f(x_{j,k})\right)
\widetilde\vartheta_{\varepsilon }(y)dy, 
\end{eqnarray*}
 
The assumption $f\in C^s(\Hn)$ implies $\nabla_\Hn^\alpha f \in C^\sigma(\Hn)$, thus ultimately providing:
$$|d_{j,k}^\varepsilon| \leq 2^{-j({[s]}+\sigma)} \norme[C^s(\Hn)]{f} \sum_{|\alpha|={[s]}}|l_\alpha|
\int_{\Hn} \norme[\Hn]{y}^\sigma |\widetilde \vartheta_{\varepsilon}(y)| dy = C 2^{-js}.$$

When ${[s]}=2m-1$ is odd ($m\geq 1$), one has:
$$
d_{j,k}^\varepsilon =
\int_{\Hn} (X\mathcal{L}^{m-1} g_{j,k}) (X\widetilde\vartheta_\varepsilon)
+\int_{\Hn} (Y\mathcal{L}^{m-1} g_{j,k}) (Y\widetilde\vartheta_\varepsilon).
$$
 Again,  $X\widetilde\vartheta_\varepsilon$ and $Y\widetilde\vartheta_\varepsilon$  are well-localized functions, with at least one vanishing moment. Using the same arguments as above,  
\begin{equation*}
d_{j,k}^\varepsilon  = \sum_{|\alpha|={[s]}} l'_\alpha 
\int_{\Hn}\left(\nabla_\Hn^\alpha g_{j,k}(y)-\nabla_\Hn^\alpha g_{j,k}(0)\right)
\widetilde\Psi_{\varepsilon,\alpha}(y)dy
\end{equation*}
with $\widetilde\Psi_{\varepsilon,\alpha} = X\widetilde\vartheta_\varepsilon$ or
$Y\widetilde\vartheta_\varepsilon$  (depending on whether the first slot in $\alpha$ codes for $X$ or $Y$) and  some  other coefficients $l'_\alpha$. The rest of the proof is the same, giving finally $|d_{j,k}^\varepsilon |\leq C2^{-js}$.


\subsection{H\"older estimate derived from wavelets coefficients}
 \label{sec32}
Let us now focus on the converse assertion in Theorem \ref{thglobal} and assume that \eqref{GLOBALHOLDER} holds.
The normal convergence of the series \eqref{WAVELETS} up to the $[s]^{\text{th}}$ derivatives
ensures that, for any multi-index $\alpha$ such that $|\alpha|\leq [s]$ the function $\nabla_{\Hn}^\alpha f$
is continuous and that the following identity holds:
\begin{equation*}
\nabla_\Hn^\alpha f(x)=\sum_{\varepsilon,j,k} 2^{j|\alpha|} d^\varepsilon_{j,k}(f)\cdot
(\nabla_\Hn^\alpha \Psi_\varepsilon)\left(2^j\circ (x_{j,k}^{-1}  x)\right).
\end{equation*}
Let us estimate the $[s]^{\text{th}}$ derivatives.
As before $s= [s]+\sigma$ and for $|\alpha|=[s]$, one gets:
$$|\nabla_\Hn^\alpha f(x)-\nabla_\Hn^\alpha f(y)| \leq
\sum_{\varepsilon,j,k} 2^{-j\sigma}
\left|
(\nabla_\Hn^\alpha \Psi_\varepsilon)\left(2^j\circ (x_{j,k}^{-1}  x)\right)
-(\nabla_\Hn^\alpha \Psi_\varepsilon)\left(2^j\circ (x_{j,k}^{-1}  y)\right)
\right|.$$
Let $j_0\in\mathbb{Z}$ such that $2^{-j_0-1} \leq \delta(x,y) < 2^{-j_0}$.
There exists $  \widetilde k=(k_1,k_2,k_3)\in\mathcal{Z}$ such that $x$ and $y$ both
belong to the dyadic neighborhood of $x_{j_0, \widetilde k}$, namely
$ x,y \in  \Lambda_{j_0,\widetilde k}$,
where $\Lambda_{j_0,\widetilde k}$ has been defined by~\eqref{DyadicNeighbor} and the remarks that follow.

\sk
 For $j\leq j_0$, one
uses that $\nabla_\Hn^\alpha \Psi_\varepsilon$ is Lipschitz:
$$
\left|
(\nabla_\Hn^\alpha \Psi_\varepsilon)\left(2^j\circ (x_{j,k}^{-1}  x)\right)
-(\nabla_\Hn^\alpha \Psi_\varepsilon)   \left(2^j\circ (x_{j,k}^{-1}  y)\right)
\right| \leq  C 2^j \delta(x,y)  \times  \xi_{j,k},  $$
where $\xi_{j,k}$ satisfies  the following estimate (Proposition~\ref{lipschitz})
\begin{align*}
\xi_{j,k} & =
\sup_{\|z\|_{\Hn}\leq \gamma_2 2^{j-j_0}} |\nabla_\Hn\nabla_\Hn^\alpha\Psi_\varepsilon ((2^j\circ (x_{j,k}^{-1}x_{j_0,\widetilde k}))\ast z) |\\
&\leq \sup_{z\in k^{-1}\ast (2^{j-j_0}\circ \mathcal{B}(\widetilde k,\gamma_2))} \left(C_0 \exp\left(- {\norme[\Hn]{z}} /{r_0}\right) \right)
 \leq C_0'  \exp\left(- {
\delta(k,2^{j-j_0}\circ \widetilde k) }/ {r_0'}\right).
\end{align*}
 
 For $j>j_0$, one uses simply the boundedness of $\nabla_\Hn^\alpha \Psi_\varepsilon$:
\begin{gather}\nonumber
\left|
(\nabla_\Hn^\alpha \Psi_\varepsilon)\left(2^j\circ (x_{j,k}^{-1}  x)\right)
-(\nabla_\Hn^\alpha \Psi_\varepsilon)\left(2^j\circ (x_{j,k}^{-1}  y)\right)
\right| \\\nonumber
\qquad\qquad\qquad\qquad \leq C_0 \left\{\exp\left(- {\delta(k, 2^j\circ x)} /{r_0}\right)+
\exp\left(- {\delta(k, 2^j\circ y)} /{r_0}\right)\right\}.
\end{gather}

Each right-hand side is obviously summable in the variable $k\in\mathcal{Z}$ because there exists a constant $C$ such that
$$
\forall z\in\Hn, \qquad \sum _{k\in \mathcal{Z}}  \exp\left(- {\delta(k,z) } / {r_0}\right) \leq C.
$$
Combining the previous estimates and the fact that the $\varepsilon$ variable belongs to a finite set $\{1,\ldots,2^Q-1\}$,
one gets the following upper bound for $|\nabla_\Hn^\alpha f(x)-\nabla_\Hn^\alpha f(y)| $:
$$
\sum_{j\leq j_0} 2^{j(1-\sigma)} \delta(x,y)+ \sum_{j>j_0}2^{-j\sigma} \leq C (2^{j_0(1-\sigma)}
\delta(x,y)+2^{-j_0\sigma})\leq C' \delta(x,y)^\sigma
$$
\textit{i.e.} $f\in C^s(\Hn)$.

 \section{Pointwise H\"older regularity}
 \label{secpointwise}

The wavelet decay property \ref{PointRegWav} of functions $f$ belonging to  $ C^s(x_0)$ is obtained in Section \ref{sec41}.  The second part of Theorem \ref{thlocal} is more delicate and is explained in Section~\ref{sec42}. Finally, Section \ref{genericcalpha} contains the proof of Theorem \ref{thgenericCalpha}, which gives the existence of a generic set of functions in $ C^s$ with a pointwise H\"older exponent everywhere equal to~$s$. This proves in turn the optimality of Theorem \ref{thglobal}.

\subsection{Upper bound for the wavelets coefficients}

\label{sec41}
 Assume first that $f\in C^s(x_0)$.
 Let $P$ be the unique polynomial of degree $\operatorname{deg}_{\Hn}(P)<s$ and $\eta>0$ such that
\eqref{Ch} holds on a small neighborhood of the origin $\mathcal{B}(0,\eta)$. 
Using the vanishing moments of  $\Psi^\varepsilon$, one has
$$ d^\varepsilon_{j,k}(f) =  2^{jQ} \int_{\Hn} (f(x_0   x) - P(x))\Psi^\varepsilon_{j,k}(x_0  x)dx$$
thus
\begin{eqnarray*}
| d^\varepsilon_{j,k}(f)| &  \leq &  C2^{jQ}  \int_{\mathcal{B}(0,\eta)}  \|x\|_{\Hn}^s  \left| \Psi^\varepsilon_{j,k}(x_0  x)\right|  dx    + 2^{jQ}  \int_{\Hn\setminus \mathcal{B}(x_0,\eta)} |f(x)|  \left| \Psi^\varepsilon_{j,k}(x)\right|  dx \\
&& + \ 2^{jQ}  \int_{\Hn\setminus \mathcal{B}(x_0,\eta)}  | P(x_0^{-1}x ) |    \left| \Psi^\varepsilon_{j,k}(x)\right|  dx.
\end{eqnarray*}
We denote the three integrals above by respectively $I_1$, $I_2$ and $I_3$.
In the following, we assume (as in the statement of Theorem~\ref{thlocal}) that \begin{equation}\label{choseR}\delta(x_{j,k},x_0)< R\end{equation}
for some constant $R\leq 1$ that we will adjust on the way. Observe that $R$ will not depend on $j$ or $k$.

\sk
The change of variables $y=2^j\circ (x_{j,k}^{-1}   x_0   x)$,  H\"older's inequality $\|a   b\|_{\Hn}^s \leq C_s (\|a\|^s_{\Hn} +\|  b\|^s_{\Hn})$ and the exponential decay~\eqref{DECAYPsi} of the mother wavelet yield:
\begin{eqnarray*}
I_1 &  = &  C  \int_{\mathcal{B}(2^j(x_{j,k}^{-1}  x_0), 2^j\eta )}  \| x_0^{-1}   x_{j,k}  (2^{-j} \circ x)\|_{\Hn}^s  \left| \Psi_\varepsilon (x )\right|  dx   \\
& \leq & C  \int_{\Hn}
 C_s\left( \delta( x_0, x_{j,k}) ^s + \|2^{-j} \circ x\| _{\Hn}^s\right) \times C_0\exp\left(   - {\|x\|_{\Hn}} /{r_0}\right)  dx   \\
& \leq & C \left( \delta( x_0, x_{j,k}) ^s +  2^{-js}   \right)   \end{eqnarray*}
for some other constant $C$, independent of $x_0$, $j$ and $k$.  The decay property \eqref{DECAYPsi} of the wavelet was used to obtain the second line.

 \sk
 For the second integral, one   uses  the Cauchy-Schwarz inequality:
 \begin{eqnarray*}
I_2 & = &    \int_{\Hn \setminus \mathcal{B}(x_0,\eta)} |f(x)|  \times 2^{jQ}\left| \Psi^\varepsilon_{j,k}(x )\right|  dx   \ \leq \  \|f\|_{L^2(\Hn)}  \left(  \int_{\Hn\setminus \mathcal{B}(x_0,\eta)} 2^{2jQ} \left| \Psi^\varepsilon_{j,k}(x )\right|^2  dx\right)^{1/2}.
 \end{eqnarray*}
The weight of the tail of the wavelet depends  on how $\delta(x_{j,k},x_0)$ compares to $\eta$.
Let us assume that $R=\eta/C$ in \eqref{choseR} with the constant $C$ given by  Corollary~\ref{CorTRIANGULAR}, thus
$$\mathcal{B}(0,2^j\eta/C)\subset \mathcal{B}(2^j\circ (x_{j,k}^{-1}  x_0),2^j \eta).$$
Observe  that  $\norme[\Hn]{y}\geq 2^j\eta/C$ implies $1\leq 2^{-j}    \norme[\Hn]{y} C/\eta$
and $  \norme[\Hn]{x_0^{-1}x_{j,k}}\leq R =\eta/C\leq 2^{-j}\norme[\Hn]{y}$.
Then, the usual change of variable $y=2^j\circ (x_{j,k}^{-1} x)$ reads
 \begin{eqnarray*}
 \int_{\Hn \setminus \mathcal{B}(x_0,\eta)} 2^{2jQ} \left| \Psi^\varepsilon_{j,k}(x )\right|^2  dx & =  & 2^{jQ}\int_{\Hn 
 \setminus \mathcal{B}(2^j\circ (x_{j,k}^{-1}  x_0),2^j \eta)}   \left| \Psi_\varepsilon(y)\right|^2  dy\\
 &  \leq &  2^{jQ}\int_{ \{y\in\Hn \,:\, \|y\|_{\Hn} \geq 2^j\eta/C\}}   \left(\frac{C \norme{y}}{2^j\eta }\right)^{Q+2s} \left| \Psi_\varepsilon(y)\right|^2  dy\\
 & \leq  & C2^{-2js}.
\end{eqnarray*}
{The last inequality uses the decay property~\eqref{DECAYPsi} of $\Psi^\ep$.} Finally, one gets $I_2\leq C 2^{-js}$.


\sk
For $I_3$, the idea is similar except that one has to compensate for $P\not\in L^2(\Hn)$ by adding an extra weight.
For example, one choses an integer $N> s$ such that $(1+\norme[\Hn]{x})^{-N}P(x)$ is bounded. Then, as previously, one has:
\begin{eqnarray*}
I_3 &\leq  & 2^{jQ}  \norme[L^\infty(\Hn)]{(1+\norme[\Hn]{\cdotp})^{-N}P}\times\int_{\Hn\setminus \mathcal{B}(x_0,\eta)}  (1+\|x_0^{-1} x\|_\Hn)^{N}   \left| \Psi^\varepsilon_{j,k}(x)\right|  dx\\
&\leq &  C \int_{\norme[\Hn]{y}\geq 2^j\eta/C}  \left(1+\|x_0^{-1} x_{j,k}\|+2^{-j}\|y\|_\Hn\right)^{N}  \left| \Psi^\varepsilon(y)\right|  dy
\end{eqnarray*}
Hence,  using the same arguments as those developed for $I_2$, the estimates boil down to:
\begin{align*}I_3 &\leq C\int_{\norme[\Hn]{y}\geq 2^j\eta/C}  \left(2^{-j}\|y\|_\Hn\right)^{N}  \left| \Psi^\varepsilon(y)\right|  dy \leq  C 2^{-jN}\int_{\Hn} \|y\|_\Hn^N  \exp\left(-\frac{\norme[\Hn]{y}}{r_0}\right) dy
\end{align*}
Finally, one gets $I_3 \leq C 2^{-jN}\leq C 2^{-js}$.
 
 \sk
Putting together the estimates for $I_1$, $I_2$ and $I_3$, one gets: 
 $$|d^\varepsilon_{j,k}(f)| \leq C \left(2^{-js}  + \delta( x_{j,k},   x_0)^{s} \right)$$
 when $ \delta( x_{j,k},   x_0)<\eta/C=R$ which is equivalent to \eqref{PointRegWav}. This concludes the proof of the 
 first half of Theorem~\ref{thlocal}.\cqfd

 \subsection{Pointwise H\"older estimate derived from wavelets coefficients}
 \label{sec42}
 Let us move to the second part of Theorem~\ref{thlocal} and prove the converse property.
One assumes that $f\in C^\sigma(\Hn) $ for some $\sigma>0$ and that \eqref{PointRegWav} holds for all triplets $(\ep,j,k)$ such that $\delta(x_{j,k},x_0)  \leq R$ for some $R>0$.  
 Let us fix $x$ such that $\delta(x,x_0)\leq R$ and let $j_0$ and $j_1$ be the unique integers such that 
 \begin{equation}
 \label{defj0j1}
 2^{-j_0-1}\leq \delta(x,x_0)<2^{-j_0} \ \text{ and } \ j_1=\left[ \frac{s}{\sigma} \cdotp j_0\right].
 \end{equation}
We aim at proving \eqref{pertelog} for $x$ close enough to $x_0$ \textit{i.e.} for $j_0$ large enough.

The wavelet decomposition of $f$ is $\displaystyle f=\sum_{j\in\mathbb{Z}} f_j(x) $, where for every $j\in \Z$, 
$$f_j(x) =\sum_{\varepsilon\in\{1,...,15\}}  \ \sum_{k\in \mathcal{Z}}  \ d^\varepsilon_{j,k}(f)\Psi^\varepsilon_{j,k} (x).$$
For subsequent use, let us notice immediately that the low frequency term
$$
 f^{\flat}(x)=\sum_{j=-\infty}^{0} f_j(x)
$$
is as regular as the wavelet itself. In particular at least $C^{[s]+2}(\Hn)$.

\mk
Assumption   \eqref{PointRegWav} reads
\begin{equation}\label{majwav}
| d^\varepsilon_{j,k}(f)| \leq C ( 2^{-js}  + \|  x_{j,k}^{-1}   x_0\|^{s}_{\Hn}) \leq C ( 2^{-js}  + \|  x_{j,k}^{-1}   x \|^{s}_{\Hn}+ \|  x ^{-1}   x_0\|^{s}_{\Hn}) .
\end{equation}
For every $n\in \{0,\ldots, [s]+1\}$ and any multi-index $\alpha$ with $|\alpha|=n$, one has:
\begin{equation}
\label{eqfj}
\nabla_\Hn^\alpha f_j(x)  =      \sum_{\varepsilon } \sum_{k\in \mathcal{Z}} d^\varepsilon_{j,k}(f) \cdot (\nabla_\Hn^\alpha \Psi^\varepsilon_{j,k})(x).
\end{equation}
As $\Psi^\ep_{j,k}(x) = \Psi_\ep(2^j\circ (x_{j,k}^{-1}  x))$ and using \eqref{DECAYPsi}, a computation similar to~\eqref{DerivativesWithScaling} gives:
\begin{equation}
\label{eqq1}
| \nabla_\Hn^\alpha \Psi^\ep_{j,k}(x)| \leq   C 2^{j|\alpha|}  \exp\left(-  {\|2^j\circ (x_{j,k}^{-1}  x)\|_{\Hn}}/{r_0}\right)  \leq \frac{C 2^{j|\alpha|}}{(1+\|2^j\circ (x_{j,k}^{-1}  x)\|_{\Hn})^{Q+1+s}}\cdotp
\end{equation} 

\mk
Next, let us notice that for a constant that does not depends on $j\in\mathbb{Z}$ or $x\in\Hn$,
 \begin{equation}
 \label{maj3}
 \forall \gamma >0, \quad\exists C>0, \qquad
 \sum_{k\in \mathcal{Z}}    \frac{   \|  x_{j,k}^{-1}   x \|^{\gamma}_{\Hn} } {(1+\|2^j\circ (x_{j,k}^{-1}   x)\|_{\Hn})^{Q+1+\gamma}} \leq C 2^{-j\gamma}.  
 \end{equation}

Combining \eqref{eqq1},~\eqref{majwav} and~\eqref{maj3} provides, in a neighborhood of $x_0$:
\begin{equation}  \label{majfjalpha}
 |\nabla_\Hn^\alpha  f_j (x)|   \leq    C 2^{j|\alpha|} 
  \sum_{k\in \mathcal{Z}}    \frac{( 2^{-js}  + \|  x_{j,k}^{-1}     x \|^{s}_{\Hn} +  \|  x^{-1}     x_0\|^{s}_{\Hn}) }{(1+\|2^j\circ (x_{j,k}^{-1}  x)\|_{\Hn})^{Q+1+s}}  \leq   C2^{j|\alpha|} \left(2^{-js} +  \|  x ^{-1}   x_0\|^{s}_{\Hn}\right).
 \end{equation}
 In particular $ | \nabla^\alpha f_j (x_0)| \leq    C 2^{j(|\alpha|-s)} $ and the series (used subsequently) $\sum_{j=0}^{\infty}      \nabla_\Hn^\alpha f_j   (x_0)   $
converges absolutely for every $\alpha$ such that $|\alpha|\leq [s]$.

Let us now introduce the (right)-Taylor polynomial $P_j$ of $f_j$ at $x_0$. According to~\eqref{TAYLOR}, it can be written
\begin{equation}
\label{eqPj}
P_j(y)= \sum_{|\alpha|=|\beta|\leq[s]}  c_{\alpha,\beta} \nabla_\Hn^\beta  f_j (x_0) y^\alpha.
\end{equation}
The coefficients $c_{\alpha,\beta}$ are chosen once and for all for the rest of this computation.
Let also $P^\flat(y)$ stand for the (right)-Taylor polynomial of  the low frequency part $f^\flat$ at $x_0$.

The polynomial $P$ that we are going to use to prove~ \eqref{pertelog} is defined by:
 $$P(y) =P^\flat(y) + \sum_{j=0}^{\infty} P_j(y).$$
 The previous estimates ensure that $P$ is indeed well defined and of degree at most $[s]$. One gets the following decomposition:
\begin{equation}\label{CSlog} 
|f(x) -P(x_0 ^{-1} x)|  \leq  |f^\flat(x)-P^\flat(x_0^{-1}x)|+  \sum_{j=0}^{j_0}  \left|  f_j(x) - P_j(x_0^{-1}  x)  \right|  + R_1(x) + R_2(x)
 \end{equation}
with two remainders:
$$R_1(x) = \sum_{j= j_0}^{\infty} |f_{j}(x) |  \quad\text{and}\quad  R_2(x)= \sum_{j= j_0}^{\infty} \left|  P_j(x_0^{-1}  x) \right|.$$
The low frequency is instantaneously dealt with by Theorem~\ref{taylor}:
\begin{equation}\label{CSlog1}
|f^\flat(x)-P^\flat(x_0^{-1}x)|\leq C \norme[\Hn]{x_0^{-1}x}^{[s]+1}.
\end{equation}

Let us now focus on the three other terms.
One uses the Taylor development of the wavelet at $x_0$ and
the unicity of the Taylor expansion to recover the polynomial $P_j$. Let us thus write:
$$\Psi_{j,k}^\varepsilon(x)=\sum_{|\alpha|=|\beta|\leq[s]} c_{\alpha,\beta}\nabla_\Hn^\beta\Psi_{j,k}^\varepsilon(x_0)(x_0^{-1}x)^\alpha + R_{j,k}^{\varepsilon}(x).$$
Theorem~\ref{taylor} ensures that, for some constant $r_1>0$ and $x$ in the neighborhood of $x_0$:
\begin{equation}\label{talorwavelet}
|R_{j,k}^\varepsilon(x)|\leq  C \norme[\Hn]{x_0^{-1}x}^{[s]+1} \sup_{\substack{|\alpha|=[s]+1, \,  \norme{z}\leq r_1}} \left|\nabla_\Hn^\alpha \Psi_{j,k}^{\varepsilon}(x_0z)\right|.
\end{equation}
Substitution in the definition of $f_j$ reads:
$$f_j(x) = \sum_{\varepsilon,k} d_{j,k}^\varepsilon(f) \left(\sum_{|\alpha|=|\beta|\leq[s]} c_{\alpha,\beta}\nabla_\Hn^\beta\Psi_{j,k}^\varepsilon(x_0)(x_0^{-1}x)^\alpha + R_{j,k}^{\varepsilon}(x)\right).$$
In the first double sum, by combining \eqref{eqfj} and \eqref{eqPj}, one recognizes $P_j(x_0^{-1}x)$ and thus
\begin{equation*}
f_j(x)-P_j(x_0^{-1}x)= \sum_{\varepsilon,k} d_{j,k}^\varepsilon(f) R_{j,k}^\varepsilon(x).
\end{equation*}
Combining~\eqref{majwav}, \eqref{talorwavelet} and the definition~\eqref{defj0j1} of $j_0$ gives:
\begin{align*}
\left|f_j(x)-P_j(x_0^{-1}x)\right| &\leq \sum_{\varepsilon,k} |d_{j,k}^\varepsilon(f)| |R_{j,k}^\varepsilon(x)| \\
&\leq C 2^{-j_0([s]+1)} \sum_{\varepsilon,k} ( 2^{-js}  +2^{-j_0s}+ \|  x_{j,k}^{-1}   x \|^{s}_{\Hn})\sup_{\substack{|\alpha|=[s]+1\\\norme{z}\leq r_1}} \left|\nabla_\Hn^\alpha \Psi_{j,k}^{\varepsilon}(x_0z)\right|.
\end{align*}
 To deal with the summation in $k$, one uses \eqref{eqq1} and \eqref{maj3}: for all $\ep=1,..., 2^Q-1$, 
$$\sum_k \sup_{\norme{z}\leq r_1} \left|\nabla_\Hn^\alpha \Psi_{j,k}^{\varepsilon}(x_0z)\right| \leq
C\sum_k  \frac{2^{j|\alpha|}}{(1+\|2^j\circ (x_{j,k}^{-1}  x_0)\|_{\Hn})^{Q+1+s}} \leq C 2^{j|\alpha|}$$
and similarly
$$\sum_k  \|  x_{j,k}^{-1}   x \|^{s}_{\Hn}\sup_{\norme{z}\leq r_1} \left|\nabla_\Hn^\alpha \Psi_{j,k}^{\varepsilon}(x_0z)\right| \leq
C\sum_k  \frac{2^{j|\alpha|} \|  x_{j,k}^{-1}   x_0 \|^{s}_{\Hn}}{(1+\|2^j\circ (x_{j,k}^{-1}  x_0)\|_{\Hn})^{Q+1+s}}\leq C2^{j(|\alpha|-s)}.$$
The summation in $\ep$ plays no role.  Putting it all together, one gets:
\begin{equation*}
\left|f_j(x)-P_j(x_0^{-1}x)\right|  \leq C 2^{-(j_0-j)([s]+1)} \left( 2^{-js}  +2^{-j_0s} \right).
\end{equation*}
Finally, the sum over $j\in\{0,\ldots,j_0\}$ boils down to:
\begin{align}
\notag
\sum_{j=0}^{j_0}\left|f_j(x)-P_j(x_0^{-1}x)\right|  &\leq C 2^{-j_0([s]+1)}\sum_{j=0}^{j_0}2^{j([s]+1-s)} + 2^{-j_0([s]+1+s)}\sum_{j=0}^{j_0}2^{j([s]+1)}\\
&\leq C 2^{-j_0s} \leq C \norme[\Hn]{x_0^{-1}x}^s.\label{CSlog2}
\end{align}

%
%
%

The term  $R_1$ contains the high-frequency components of the Littlewood-Paley decomposition of~$f$ and
is  responsible for the logarithmic correction in~\eqref{pertelog}.
By  \eqref{defj0j1} and~\eqref{majfjalpha}, 
$$\forall \ j\geq j_0,  \ \  |f_j(x)| \leq C(2^{-js }+  \|  x ^{-1}   x_0\|^{s}_{\Hn}) \leq C (2^{-js}+2^{-j_0s})\leq  C   \|  x ^{-1}   x_0\|^{s}_{\Hn}.$$
Let us split this remainder depending on whether $j_0 \leq j < j_1$ or $j\geq j_1 $
$$R_1(x) \leq  \sum_{j = j_0}^{j_1} |f_{j}(x) | +\sum_{j = j_1}^{\infty} |f_{j}(x) |.$$

Our choice  \eqref{defj0j1} for $j_1$ and $j_0$  gives $j_1 \sim  s j_0/\sigma \sim s/\sigma \cdot |\log \|  x ^{-1}   x_0\|^{s}_{\Hn}|$. Hence
\begin{eqnarray}
\label{final1}  \sum_{j = j_0}^{j_1} |f_{j}(x) |  & \leq &     j_1 \cdot C   \|  x ^{-1}   x_0\|^{s}_{\Hn} \leq C  \|  x ^{-1}   x_0\|^{s}_{\Hn}  \cdot |\log  \|  x ^{-1}   x_0\|^{s}_{\Hn}|.
 \end{eqnarray}
When $j\geq j_1$, one uses $f\in C^\sigma (\Hn)$ instead and  Theorem~\ref{thglobal} which gives  $|d_{j,k}^\ep(f)|\leq C 2^{-j\sigma}.$
Combined with~\eqref{DECAYPsi}, one deduces that 
$$ |   f_j (x)|    \leq   \sum_{\ep}   \sum_{k\in \mathcal{Z}} C 2^{-j\sigma}  e^{-\delta(x_{j,k},x) / r_0} \leq C 2^{-j\sigma}.$$
Using  \eqref{defj0j1} one last time yields to the expected conclusion:
\begin{equation}
\label{final2}
\sum_{j = j_1}^{\infty} |f_{j}(x) |   \leq  C\sum_{j =  j_1}^{\infty} 2^{-j\sigma} \leq C 2^{-j_1\sigma} \leq C \norme[\Hn]{x_0^{-1}x}^s.
\end{equation}

Let us move to $R_2$ which contains the Taylor expansions of the high-frequency components of the Littlewood-Paley decomposition of~$f$.
Intuitively, it is small because of the natural spectral separation between polynomials and highly-oscillatory functions.

Using \eqref{defj0j1} and \eqref{majfjalpha}, each term of the sum boils down to:
$$\left|  P_j(x_0^{-1}  x) \right| \leq \sum_{|\alpha|=|\beta|\leq[s]} |c_{\alpha,\beta}| |\nabla_\Hn^\beta f_j(x_0)| \|(x_0^{-1}x)^\alpha\|_\Hn
\leq C \sum_{n=0}^{[s]} 2^{j(n-s)-j_0n}$$
and thus
\begin{align} 
\sum_{j= j_0}^{\infty} \left|  P_j(x_0^{-1}  x) \right| &\leq C \sum_{n=0}^{[s]}2^{-j_0n}\left(\sum_{j\geq j_0} 2^{-j(s-n)}\right)  \leq C 2^{-j_0s} \leq C \norme[\Hn]{x_0^{-1}x}^s.
\label{finpointwise2}
\end{align}

Substituting   \eqref{CSlog1}, \eqref{CSlog2}, \eqref{final1},  \eqref{final2} and \eqref{finpointwise2} back in the original
question \eqref{CSlog} proves that  \eqref{pertelog} holds in a neighborhood of $x_0$ and concludes the proof of
Theorem~\ref{thlocal}. \cqfd

 \subsection{Generic monofractactality of functions in $C^s(\Hn)$}
\label{genericcalpha}

The proof of Theorem~\ref{thgenericCalpha} is classical in the Euclidian context \cite {Jaffard_generic} and can be adapted quickly to ours. 
\mk

Let us recall that for any   $f\in C^s(\Hn)$,  Theorem \ref{thglobal} gives a constant $C>0$ such that
\begin{equation}\label{WAVELETS2}
f=\sum_{\varepsilon,j,k} d^\varepsilon_{j,k}(f)\Psi^\varepsilon_{j,k}
\qquad\text{with}\qquad
| d^\varepsilon_{j,k}(f) | \leq C 2^{-js}
\end{equation}
and  $\|f\|_{C^s} = \inf\{C>0: \eqref{WAVELETS2} \mbox{ is satisfied for all }\ep, \, j, \, k\}$
is a Banach norm on $C^s(\Hn)$. 
For each integer $N$, let us define:
\begin{equation}\label{defFN} 
\begin{aligned}
E_N  &  =  \left\{f\in C^s(\Hn) \, : \, \forall \,(\ep, \, j, \, k), \enspace 2^{js+N} d_{j,k}^\ep(f) \in \mathbb{Z}^\ast \right\}\\
 F_N  & =   \left\{g\in C^s(\Hn) \, : \, \exists f\in E_N, \enspace \norme[C^s(\Hn)]{f-g} < 2^{-N-2}\right\}.
\end{aligned}
\end{equation}

\begin{lemma}
For every $N\geq 1$, all functions in $F_N$ are monofractal of exponent $s$: $$\forall g\in F_N, \qquad \forall x\in\Hn,\qquad h_g(x)=s.$$
\end{lemma}
\Proof
This simply follows from the fact that, given $f\in E_N$, all the wavelet coefficients of $f$ satisfy
$$ 2^{-N-js} \leq |d_{j,k}^\ep(f)| \leq \|f\|_{C^s}2^{-js}  .$$
Thus for any function $g\in F_N$ and its associated $f\in E_N$:
$$ 2^{-N-js} - 2^{-N-2-js}   \leq |d_{j,k}^\ep(g)| \leq  \|f\|_{C^s}  2^{-js} +2^{-N-2-js}$$
\textit{i.e.}
$$ 2^{-N-1-js}   \leq |d_{j,k}^\ep(g)| \leq  \left(\|f\|_{C^s} +2^{-N-2}\right) 2^{-js}.$$
In particular, $g\in C^s(x)$ for any $x\in \Hn$ and  there is no $x_0\in\Hn$ and $s'>s$ such that $g\in C^{s'}(x_0)$. Indeed,    \eqref{PointRegWav} with $s'>s$ is not compatible when $j$ tends to infinity with the left hand-side of the above inequality. 
\cqfd

\begin{lemma}
The set  $
 \mathcal{R}   =  \bigcup_{N\geq 1} \ F_N $
 is a dense open set in $C^s(\Hn)$ containing only monofractal functions with exponent $s$.
\end{lemma}
\Proof
 The preceding lemma ensures that $\mathcal{R}$ is composed of monofractal functions.
 According to  \eqref{defFN}, $F_N$ is  an open set and thus, so is $\mathcal{R}$. 
Let us check the density.
Given $f  \in C^s(\Hn)$ and $\eta>0$, let us choose $N \in\mathbb{N}$ so that $2^{-N} < \eta$.
Let us define the ``non-zero integer part'' function
$$E^\ast(x)=
\begin{cases}
1&\text{if }0\leq x <2,\\
 [x]  & \text{else}.
\end{cases}$$
Obviously $E^\ast : \R\to\Z^\ast$ and $|x-E^\ast(x)|\leq 1$.
Let us finally define a function $g\in F_N$ by its wavelets coefficients:
$$d_{j,k}^\varepsilon(g) = 2^{-js-N}E^\ast\left( 2^{js+N}d_{j,k}^\varepsilon(f) \right).$$
By construction, 
$$2^{js}\left|d_{j,k}^\varepsilon(f)-d_{j,k}^\varepsilon(g)\right| = 2^{-N}|2^{js+N}d_{j,k}^\varepsilon(f) -E^\ast\left( 2^{js+N}d_{j,k}^\varepsilon(f) \right)|\leq 2^{-N}<\eta$$
thus $\|f-g\|_{C^s} < \eta$. This proves the density of  $\mathcal{R} $ in $C^s(\Hn)$.
 \cqfd

\section{Upper bound  for the  multifractal spectrum   in a Besov space}
\label{secbesov1}

The classical Sobolev embedding $B^{s}_{p,q}(\Hn) \hookrightarrow C^{s-Q/p}(\Hn)$ can be retrieved easily using wavelets. 
Indeed, the definition \eqref{defbesov} reads
$$\norme[\ell^p(k)]{2^{j\left(s- Q/p\right)}d^\varepsilon_{j,k}(f)} \in \ell^q(j)$$
and implies the existence of a constant $C_0>0$ such that for every triplet $(\ep,j,k)$:
\begin{equation}\label{sobolev_embedding}
|d_{j,k}^\ep(f)|\leq C_0 2^{-j(s-Q/p)}.
\end{equation}
Thus~\eqref{GLOBALHOLDER} holds and Theorem~\ref{thglobal} ensures that $f\in C^{s-Q/p}(\Hn)$. In particular, \eqref{PointRegWav} is satisfied around any point $x_0\in \Hn$.
and thus by Theorem~\ref{thlocal}, one has
$$\forall x\in\Hn, \qquad h_f(x)\geq s-Q/p.$$
It is worth mentioning that the index $q$ of the Besov space $B^s_{p,q}(\Hn)$ does not play any role in the Sobolev embedding and neither does it in the multifractal analysis of $f$.

\mk

Let us now establish Theorem~\ref{majspectrumBesov}, \textit{i.e.} that for any $h\geq s- Q/p$, the iso-H\"older set of regularity $h$ is of Hausdorff dimension
$$d_f(h)\leq \min\left(Q \,,\, p\left(h-s+ Q/p\right)\right).$$
The inequality is obvious as soon as $h\geq s$, since the upper bound reduces to $Q$ which is the Hausdorff dimension of $\Hn$ itself.
Thus one can now assume that  \begin{equation}\label{h_range}s-Q/p \leq h<s.\end{equation}
and in particular, that $1\leq p<+\infty$.

\mk

By Theorem  \ref{thlocal}, heuristically,  the   wavelet coefficients   that might give rise  to an exponent $h_f(x_0)\leq  h$ for some $x_0\in\Hn$ satisfy $|d_{j,k}^\ep| \geq 2^{-jh}$. Hence we focus on
$$N_f(j,h) =\Big \{k\in \mathcal{Z} \,:\,  \exists \  \ep\in \{1,2, ...,2^Q-1\}, \quad |d_{j,k}^\ep(f)| \geq C_0 2^{-jh} \Big\}.$$
We focus on the set $N_f(j,h) $ obtained by taking  the constant $C_0$ to be the one from the Sobolev embedding~\eqref{sobolev_embedding}. This choice is made for technical reasons that will be clear later on (see equation \eqref{calc12}).

\begin{lemma}
\label{lemma1}
There exists $C>0$ such that for every $j\geq 1$, for every $h'\in (s-Q/p,s]$,  
\begin{equation*}
 \#N_f(j,h') \leq C 2^{jp(h'-s+Q/p)}.
 \end{equation*}
\end{lemma}

\Proof
Obviously from \eqref{defbesov}, $B^s_{p,q}(\Hn)\subset B^s_{p,\infty}(\Hn)$ and thus there is a constant $C>0$ such that  for any $j\in\mathbb{Z}$, one has
$ \sum_{k\in \mathcal{Z}} 2^{j(ps-Q)}|d_{j,k}^\ep(f)|^p   \leq  C.$
Hence, 
$$\begin{aligned}
C 2^{-j(ps-Q)}  \geq   \sum_\ep\sum_{k\in \mathcal{Z}} |d_{j,k}^\ep(f) |^p   \geq   \sum_{k\in N_f(j,h') } \sum_\ep\left| d_{j,k}^\ep(f) \right |^p  \geq C_0^p\times ( \#N_f(j,h')) \times 2^{-jph'},
\end{aligned}$$
which yields the result. Observe that it also holds when $h'\geq s$ but it is useless for our purpose.
\cqfd
 
\begin{lemma}
\label{leminf}
The set
\begin{equation*}
E^{\leq}_f(h)= \{x\in \Hn: h_f(x) \leq h\}.
\end{equation*}
has the following property:
\begin{equation} \label{maj12}
\dim_H \, E^{\leq}_f(h) \leq p\left(h-s+ Q/p\right).
\end{equation}
\end{lemma}
Estimate~\eqref{maj12} is stronger than \eqref{BesovSpectrum} since  $E _f(h)  \subset E^{\leq}_f(h)  = \bigcup_{h'\leq h} E _f(h')$. 
In particu\-lar, $\dim_H E_f(h) \leq  \dim_H E^{\leq}_{f}(h)  $ and Theorem \ref{majspectrumBesov} follows immediately.

\Proof
The definition~\eqref{def_hf} of $h_f$ as a supremum implies that
$$E^\leq_f(h) = \left\{ x\in\Hn\,:\, \forall \ h'>h,\enspace f\not\in C^{h'}(x) \right\}.$$

Joint with Theorem~\ref{thlocal}, this observation provides:
\begin{equation}\label{Eleq_wavelet}
E^\leq_f(h) = \left\{ x\in\Hn\,: \,  \forall \ h'>h, \enspace \sup_{\varepsilon,j,k} \left[2^{-j}+\delta(x,x_{j,k})\right]^{-h'} |d_{j,k}^\varepsilon(f)| =+\infty\right\}.
\end{equation}
Note that as soon as $f\in L^2(\Hn)$, one has $|d_{j,k}^\varepsilon|\leq C 2^{jQ/2}$ and thus the only real constraint contained in~\eqref{Eleq_wavelet} concerns the regime $j\to+\infty$ and $h'\in[h,\alpha h]$ for an arbitrary $\alpha>1$.

\mk
We fix $\alpha>1$, and let $T$ be a  large integer. As said above, the only interesting wavelet coefficients are those satisfying $|d^\ep_{j,k}| \geq C_0 2^{-j h'}$ (it is not enough to consider only  those greater than $ C_0 2^{-j h}$). One splits $[s-Q/p,\alpha h]$ into intervals of length 
$$\eta=  \left(\alpha h-s+ Q/p\right)/T,$$
namely the intervals $I_m=[h_{m-1},h_m]$ for $m \in\{1,..., T\}$ and  $h_m=s-Q/p+m\eta$.
One choses $T$ large enough so that $h_{0}-\eta = s-Q/p-\eta>0$.
The next idea is that if $x\in\Hn$ is far   from the dyadic set $(x_{j,k})_{j\in\Z,k\in\mathcal{Z}}$ and simultaneously the wavelet coefficient  $|d^\ep_{j,k}|$ is too small, then it cannot contribute to~\eqref{Eleq_wavelet}. More precisely, if 
\begin{equation}\label{calc12}
|d_{j,k}^\ep(f)|\leq C_0 \, 2^{-jh\alpha}
\quad\text{or}\quad
\exists \ m\in \{0,\ldots,n_1-1\}, \quad \begin{cases} |d_{j,k}^\ep(f)| \leq C_0 2^{-jh_{m}} \\[.5ex] \delta(x,x_{j,k})\geq 2^{-j(h_{m}-\eta)/\alpha h}\end{cases}
\end{equation}
with the same constant $C_0>0$ as in the Sobolev embedding~\eqref{sobolev_embedding}, then
$$\left[2^{-j}+\delta(x,x_{j,k})\right]^{-h'} |d_{j,k}^\varepsilon(f)|  \leq \begin{cases}
C_0 \, 2^{-j(\alpha h-h')} & \text{in the first case,}\\
C_0 \, 2^{j(h_{m}-\eta)(h'/\alpha h)}2^{-jh_{m}} & \text{in the others}.
\end{cases}$$
In the range of $1<h'/h\leq \alpha$ and for large $j$, one infers in both cases that:
$$\left[2^{-j}+\delta(x,x_{j,k})\right]^{-h'} |d_{j,k}^\varepsilon(f)|  \leq C_0.$$
Therefore, if~\eqref{calc12} holds for \textit{any} family $(\varepsilon_n,j_n,k_n)$ with $j_n\to+\infty$, then  $x\notin E^\leq_f(h)$.

By contraposition,  for any $x\in E^\leq_f(h)$, there exists a family $(\varepsilon_n,j_n,k_n)$ with $j_n\to+\infty$  contradicting ~\eqref{calc12}.
By the Sobolev embedding~\eqref{sobolev_embedding}, each wavelet coefficient is bounded above by $C_0 2^{-jh_0}$.
Hence,  for each $n$, there exists necessarily  $m\in\{1,\ldots,T\}$ such that
$$ C_0 2^{-j_nh_{m}} < |d_{j_n,k_n}^{\ep_n}(f)|\leq C_0 2^{-j_nh_{m-1}} \quad\text{and}\quad x\in\mathcal{B}(x_{j_n,k_n},2^{-j_n(h_{m-1}-\eta)/\alpha h}).$$

The previous statement can be expressed more easily in term of lim-sup sets:
\begin{equation}
\label{majeleq}
\begin{gathered}
E^{\leq}_{f}(h) \  \subset \  \bigcup_{m=1}^{T} \ S_{m,\eta},  \ \text{ with } \ 
S_{m,\eta}\enspace=\enspace\bigcap_{J\geq 1}  \ \bigcup_{\substack{j \geq J, \ k\in N_f(j,h_m)}} \mathcal{B}(x_{j,k},2^{-j( h_m-2\eta)/\alpha h}).
\end{gathered}
\end{equation}

Let us now establish an upper bound for  the Hausdorff dimension of each set $S_{m,\eta}$. 
Given $\xi>0$, one chooses an integer $J_\xi$ so large that $2\gamma_1 \times 2^{-J_\xi (s-Q/p-\eta)/\alpha h} \leq \xi$ (with the constant $\gamma_1$ from~\eqref{TRIANGULAR}).
A covering of $S_{m,\eta} $ by balls of diameter less than $\xi$ is thus provided by:
$$S_{m,\eta} \subset \ \ \bigcup_{\substack{j \geq J_\xi, \   k\in N_f(j,h_m)}} \mathcal{B}(x_{j,k},2^{-j( h_m-2\eta)/\alpha h}).$$
For any $d\geq0$, the $\mathcal{H}^d_\xi$-premeasure of $S_{m,\eta}$ can then be estimated easily:
$$\mathcal{H}^d_\xi(S_{m,\eta}) \leq  \sum _{j \geq J_\xi} \   \ \sum_{k\in N_f(j,h_m)} \left(C 2^{-j(h_m-2\eta)/\alpha h}\right)^d.$$
Using Lemma~\ref{lemma1} to estimate $ \#N_f(j,h_m) \leq C 2^{jp(h_m-s+Q/p)}$ gives:
$$\mathcal{H}^d_\xi(S_{m,\eta}) \leq C  \sum _{j \geq J_\xi} 2^{j\left[p(h_m-s+Q/p)-d(h_m-2\eta)/\alpha h\right]}.$$
This series converges when
\begin{equation} \label{majdim}
d> \frac{\alpha ph}{h_m-2\eta}\left( h_m-s+ Q/p\right)
\end{equation}
and in that case
$$\mathcal{H}^d_\xi(S_{m,\eta}) \leq C 2^{-J_\xi\left[d(h_m-2\eta)/\alpha h-p(h_m-s+Q/p)\right]}.$$
As $\xi$ tends to zero, the $d$-Hausdorff measure of $S_{m,\eta}$ is  0, which in turn implies that $\dim_H \, (S_{m,\eta}) \leq d$. Finally, optimizing for any $d$ that satisfies~\eqref{majdim} provides 
$$\dim_H \, (S_{m,\eta}) \leq \frac{\alpha ph}{h_m-2\eta}\left( h_m-s+ Q/p\right) =  {\alpha ph}\left(1-\frac{ s-Q/p-2\eta}{h_m-2\eta} \right).$$
Looking back at \eqref{majeleq} one deduces that
\begin{eqnarray*}
 \dim_H  E^{\leq}_{f}(h)  &  \leq &\max_{m=1,...,T} {\alpha ph}\left(1  -\frac{ s-Q/p-2\eta}{h_m-2\eta} \right)  \leq   \alpha ph\left(1- \frac{s-Q/p-2\eta}{\alpha h-2\eta}\right) .
 \end{eqnarray*}
The limit $\eta\to0$ provides
$$\dim_H E^{\leq}_{f}(h)   \leq \alpha ph\left(1- \frac{s-Q/p }{\alpha h }\right)  =  p\left(\alpha h - s+ Q/p\right).$$
Finally, letting $\alpha\to 1$ gives~\eqref{maj12} and Theorem~\ref{majspectrumBesov}.
 \cqfd 


 \section{ Almost all functions in $B^s_{p,q}(\Hn)$ are multifractal}
 \label{section:generic_besov}

To prove Theorem \ref{th:besovps}, one will  explicitly construct  a $G_\delta$ set of functions in $B^s_{p,q}(\Hn)$  that satisfies \eqref{BesovSpectrumequality}.
The proof is   adapted from  the one of \cite {Jaffard_generic},    modifications are due to the metric  on $\Hn$.
We   first construct a subset $\mathcal{R}_0$ of $B^s_{p,q}(\Hn)$ whose restriction to $[0,1]^3$ is generic in $B^s_{p,q}([0,1]^3)$ and satisfies~\eqref{BesovSpectrumequality}. Next, we define: $$\forall k\in\mathcal{Z}, \qquad \mathcal{R}_{k}=\{f(k^{-1}x)\,:\, f\in\mathcal{R}_0\}.$$
Finally, the intersection $\displaystyle \mathcal{R}=\bigcap_{k\in\mathcal{Z}} \mathcal{R}_k$ is generic in $B^s_{p,q}(\Hn)$ because it is a countable intersection of $G_\delta$ sets and thus still a $G_\delta$ set.
By construction, it  will still satisfy~\eqref{BesovSpectrumequality}.

 \sk

The actual proof of Theorem \ref{th:besovps} is contained in \S\ref{step3}. To build up for it, one  needs to recall a few classical definitions and results on dyadic approximation in \S\ref{step1}.
One then constructs a single function that satisfies~\eqref{BesovSpectrumequality} in \S\ref{step2}, which is the starting point for growing the set $\mathcal{R}_0$ in \S\ref{step3}.

 \subsection{Dyadic approximation in $\Hn$}
 \label{step1}

For any $j\in\mathbb{N}$, one considers the subset of indices
\begin{equation}
\label{defL0}
\mathscr{L}_0(j) = \left\{ k\in\mathcal{Z} \,:\, x_{j,k}=2^{-j}\circ k \in[0,1)^3\right\}.
\end{equation}
For later use, let us observe immediately that
\begin{equation}\label{card_L0j}
\#(\mathscr{L}_0(j)) = 2^{Qj}.
\end{equation}

\begin{dfn}
\label{defirred}
A dyadic point $x_{J,K}$ is called irreducible if $K=(K_p,K_q,K_r)$ and  at least  one of the three fractions $2^{-J}K_p$, $2^{-J} K_q$ or $2^{-2J} K_r$ is irreducible.
A point $x_{J,K}$ is called the {\em irreducible version} of $x_{j,k}$ if $x_{J,K}$ is irreducible and $x_{j,k}=x_{J,K}.$
\end{dfn}
One can check immediately that for a given  couple $(j,k)$ the corresponding irreducible couple  $(J,K)$ is unique.
Note that one may have $(j,k)=(J,K)$ but that one always has $J \leq j$.
Conversely, given an irreducible $x_{J,K}$ and $j\geq J$, there exists a unique $k\in\mathcal{Z}$ such that $x_{j,k}=x_{J,K}$, namely $k=2^{j-J}\circ K$.

Given an integer $J\in\mathbb{N}$, the number of irreducible elements $x_{J,K}\in[0,1)^3$ is:
\begin{equation}\label{card_irreducible_antecedents}
\#\left\{K\in \mathscr{L}_0(J) \, : \, x_{J,K}\text{ irreducible}\right\} = (2^Q-1) \times 2^{Q(J-1)}.
\end{equation}
Indeed, $K=(K_1,K_2,K_3)$ provides a non irreducible $x_{J,K}\in[0,1)^3$ if and only if  $0\leq K_1,K_2<2^J$, $ 0\leq K_3<2^{2J}$ and
$$K_1\equiv 0 \enspace [\text{mod }2]\quad\text{and}\quad K_2\equiv 0 \enspace [\text{mod }2] \quad\text{and}\quad K_3\equiv 0 \enspace [\text{mod }4].$$
The complementary set in $\mathscr{L}_0(J)$ is thus of cardinal $2^{J-1}\times 2^{J-1}\times 2^{2J-2} = 2^{Q(J-1)}$ 
and~\eqref{card_irreducible_antecedents} follows from~\eqref{card_L0j}.

\begin{figure}[]
\includegraphics[width=300pt,height=140pt]{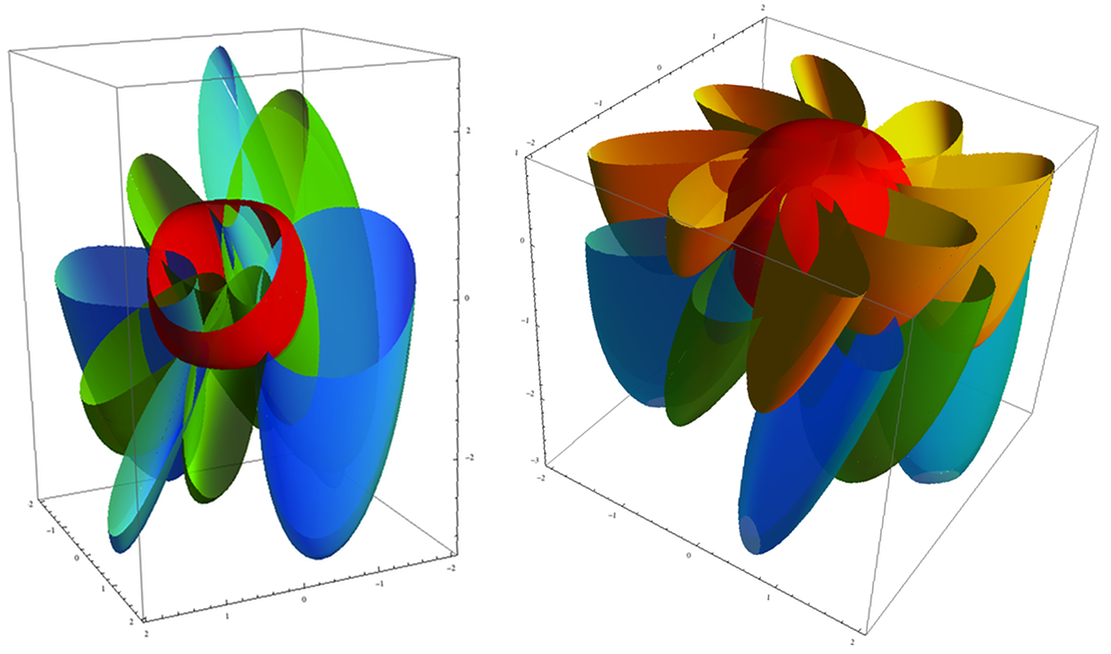}
\caption{\label{fig:balls}Left : section of the unit ball centered at the origin (red) and the balls of radius 1 centered at $(k_1,k_2,0)$ with $k_1^2+k_2^2=1$ (green) or $k_1^2+k_2^2=2$ (blue). Right : another section of the same balls (same colors) and the balls centered in $(k_1,k_2,2)$ (orange).}
\end{figure}

Recall that $\mathcal{B}(x,r)=\{y\in\Hn \,:\, \|x^{-1}y\|_{\Hn}<r\}$ denotes the open gauge ball of radius~$r$.
For fixed  $j\in\mathbb{Z}$, the dyadic elements $\{x_{j,k}: k \in \mathcal{Z}\}$ are well-distributed in $\Hn$, in the sense that the open balls $\{\mathcal{B}(x_{j,k},2^{-j}): k \in \mathcal{Z}\}$ do not intersect too much.
One can check easily the following lemma (see Figure~\ref{fig:balls}).

\begin{lemma}
\label{number cover}
For a given $j\in\mathbb{N}$ and $k \in \mathcal{Z}$, the only parameters $k'\in\mathcal{Z}$ such that $\mathcal{B}(x_{j,k k'},2^{-j}) \cap \mathcal{B}(x_{j,k},2^{-j}) \neq\emptyset$
are the 43 cubes defined by $k'=(k_1',k_2',k_3')$,  with
$$ \begin{cases}k_1'=k_2'=0\\|k_3'|\leq 1\end{cases} 
 \text{ or }  \ \ \ 
\begin{cases}1\leq {k_1'}^2+{k_2'}^2\leq 2\\|k_3'|\leq 2.\end{cases}.$$
\end{lemma}
\Proof
By scale invariance of the pseudo-norm $\|\cdot\|_{\Hn}$, it is sufficient to investigate $j=0$ and $k=(0,0,0)$. Then a counting argument applies.
\cqfd

\medskip
Observe that  if $r_0^2=\frac{\sqrt{3}}{2}<1$, then  the cylinder $\Gamma_0 = \left\{(p,q,r)\in\Hn \,:\, p^2+q^2<r_0^2\right\}$ is included in $\bigcup_{k_3\in\mathbb{Z}} \mathcal{B}((0,0,k_3), 1)$. Since for any $k=(k_1,k_2,k_3)\in\mathcal{Z}$, the left translation of $\Gamma_0$ is another vertical cylinder, 
one can choose a constant $C>0$ such that for each $j\in\mathbb{N}$, the family of balls  $$\{\mathcal{B}(x_{j,k},C2^{-j}): k \in \mathcal{Z}\}$$ covers the whole space $\Hn$
and that  for any   strictly increasing sequence $(j_m)_{m\geq 1}\in\mathbb{N}^{\mathbb{N}^\ast}$:
\begin{equation}
\label{cover}
[0,1)^3= \limsup_{m\to +\infty} \ \ \bigcup_{k\in\mathcal{L}_0(j_m)} \mathcal{B}(x_{j_m,k},C2^{-j_m}).
\end{equation}

Since each point $x\in[0,1)^3$ belongs to an infinite number of balls $\mathcal{B}(x_{j_m,k},C2^{-j_m})$, one can wonder of the exact proximity of $x$ to the dyadic elements of $\Hn$.
\begin{dfn}
\label{deftaux}
Let $\mathcal{J}=(j_m)_{m\geq 1}$ be a strictly increasing sequence of integers. For $\xi >0$, an element $x\in \Hn$ is said to be $\xi$-approximable with respect to $\mathcal{J}$ when the inequality
\begin{equation}
\label{ineg2}
\delta(x,x_{j_m,k}) \leq C 2^{- j_m \xi}
\end{equation}
holds true for an infinite number of couples $(m,k)$ and the same constant $C$ that appears in~\eqref{cover}. For a given $\xi>0$, one defines also:
$$S_{\xi}(\mathcal{J})  = \{ x \in \zu^3:    x \mbox{ is $\xi$-approximable with respect to $\mathcal{J}$ } \}.$$
The dyadic approximation rate of $x$  with respect to $\mathcal{J}$ is the real number:
$$ \xi_x( \mathcal{J}) = \sup\{\xi >0: \ x \mbox{ is $\xi$-approximable with respect to $\mathcal{J}$} \}.$$
Finally, the iso-approximable set of rate $\xi$ is:
$$\widetilde S_{\xi} ( \mathcal{J})= \{ x \in \zu^3:    \xi_x( \mathcal{J})=\xi \}.$$
\end{dfn}
When $\mathcal{J} = \mathbb{N}$ one simply writes $S_{\xi}$, $\widetilde S_{\xi}$  and $\xi_x$ and in that case it is sufficient to restrict oneself in~\eqref{ineg2} to irreducible dyadic elements $x_{j,k}$.
Let us observe that, because of \eqref{cover}:
\begin{equation*}
\forall x\in\Hn, \qquad \xi_x ( \mathcal{J})\geq 1.
\end{equation*}

\mk

The size of the sets $  S_{\xi}( \mathcal{J}) $ and $\widetilde S_{\xi}( \mathcal{J}) $ in terms of Hausdorff dimension and measures can be described thanks to the so-called {\em mass transference principle} by V.~Beresnevich, D.~Dickinson and S.~Velani \cite{BDV}.

\begin{prop}\label{propequaldim}
For every strictly increasing sequence of integers $\mathcal{J}=(j_m)_{m\geq 1}$ and every  $\xi \geq 1$,   one has:
\begin{equation}\label{equaldim}
\dim_{H}  \widetilde S_\xi ( \mathcal{J})  = \dim_H   S_\xi( \mathcal{J})    =  Q/\xi.
\end{equation}
\end{prop}

\Proof
It is quite easy to obtain that for every $\xi\geq 1$, 
\begin{equation}
\label{majdimSdelta}
\dim_H   S_\xi ( \mathcal{J})   \leq Q/\xi 
\qquad\text{and}\qquad
\dim_{H}  \widetilde S_\xi( \mathcal{J})  \leq Q/\xi
\end{equation}
Indeed, by definition, one has: $$S_\xi(\mathcal{J}) \subset \bigcap_{n \geq1} \ \ \bigcup_{\substack{j\geq n , \ k\in\mathscr{L}_0(j)}} \mathcal{B}(x_{j,k},2^{-j\xi}).$$
For $d>Q/\xi$ and an arbitrary $\eta>0$, one choses $n$ large enough such that $2^{-n\xi}<\eta$. The previous inclusion provides a covering of $S_\xi(\mathcal{J})$ by balls of radius
smaller than $\eta$, thus:
$$\mathcal{H}^d_\eta(S_\xi(\mathcal{J}))\leq C \sum_{j\geq n} 2^{jQ}(2^{-j\xi})^d \leq C 2^{n (Q-d\xi ) }\underset{n\to +\infty}{\longrightarrow}0,$$
which proves the first half. The second half follows by noticing that, as $S_{\xi'}(\mathcal{J})$ is a decreasing family (for inclusion) when $\xi'$ increases:
$$\widetilde{S}_\xi(\mathcal{J}) \subset \bigcap_{\xi'<\xi} S_{\xi'}(\mathcal{J}).$$

The converse inequality to \eqref{majdimSdelta} is very difficult,  but is contained in \cite{BDV}.
Their main theorem (stated as Theorem~\ref{thm_BDV} below) holds  at a great level of generality. It holds in particular on the Heisenberg group since $\Hn$ endowed with the metric \eqref{defmetric}
 and the Haar measure $\ell=dp\,dq\,dr$ satisfies the following conditions :
\begin{enumerate}
\item[(H1)]  $\ell$ is translation-invariant,
\item[(H2)] $\ell$ has a scaling behavior \textit{i.e.}   there exists a constant $C>1$ such that:
 $$\forall x\in \Hn, \quad \forall r>0, \qquad C^{-1} r^Q   \leq \ell(\mathcal{B}(x,r)) \leq C r^Q$$
and in particular $\ell$ is doubling \textit{i.e.}:
 $$\forall x\in \Hn, \quad \forall r>0, \qquad \ell(\mathcal{B}(x,2r)) \leq  C \ell(\mathcal{B}(x,r)).$$
 \item[(H3)]  The dyadic set $\{x_{j,k} \,  :\, j\in\mathcal{J}, \enspace k\in\mathcal{Z}\}$ is discrete.
\end{enumerate}
Joint with the covering property \eqref{cover}, the main result of \cite{BDV}, called {\em mass transference principle}, implies   the following:
\begin{thm}[\cite{BDV}, Theorem 2, p. 15]\label{thm_BDV}
For every strictly increasing sequence of integers $\mathcal{J}=(j_m)_{m\geq 1}$ and every  $\xi \geq 1$,  one has:
\begin{equation}
\label{inegmeasure}
\mathcal{H} ^{Q/\xi}(S_\xi( \mathcal{J}) )=+\infty \qquad\text{and}\qquad \mathcal{H} ^{Q/\xi}(\widetilde S_\xi( \mathcal{J}) )=+\infty.
\end{equation}
\end{thm}

The statement \eqref{inegmeasure} implies that both the Hausdorff dimension of~$S_\xi( \mathcal{J}) $  and  of $ \widetilde S_\xi ( \mathcal{J})$
are greater or equal to $Q/\xi$. Combined with~\eqref{majdimSdelta},  this proves~\eqref{equaldim}.
\cqfd

 \subsection{Example of a function with the maximal possible  spectrum}
\label{step2}

Recall that for  $x_{j,k}\in [0,1)^3$,  we denote by $x_{J,K}$ its  irreducible version.

\begin{prop}
\label{propF}
Let $\beta= {1}/{p}+ {2}/{q}$ and $F$ be the function whose wavelet coefficients are 
\begin{equation}
\label{deffjk}
F^\ep_{j,k} :=  d^\ep_{j,k}(F)=
\begin{cases}\dfrac{2^{ -j(s- Q/p )- J Q/p }}{j^\beta} & \mbox{if $x_{j,k}\in \zu^3$}\\
0 & \mbox{otherwise.}
\end{cases}
\end{equation}
The function $F $ belongs to $B^s_{p,q}(\Hn)$ and  it  satisfies \eqref{BesovSpectrumequality}.
\end{prop}
Observe that, by construction, $F$ is essentially supported in $C_0=[0,1]^3$. Outside $C_0$, $F$ is  as smooth as the mother wavelet and decays rapidly at infinity.

\mk
\Proof
Let us fix a generation $j\geq 1$ and consider the sequence   $a_j=\big\| 2^{j\left(s- Q/p\right)}  F^\varepsilon_{j,k} \big\|_{\ell^p(k)}.$
  For a given integer $j$, one has
\begin{align}
\nonumber
2^{-j(ps-Q)}a_j^p &=  \sum_{k\in \mathscr{L}_0(j)}  \sum_{1\leq\varepsilon<2^Q} \ |   F^\varepsilon_{j,k}|^p   =(  2^Q  -1)  \sum_{k\in \mathscr{L}_0(j)}    \ |   F^{1}_{j,k}|^p  \\
\label{eq22s}
& \leq   (2^Q  -1)\ \sum_{J=0}^j  \frac{ 2^{ -j(ps - Q )  -QJ  }}{j^{p\beta}} \: \#\{ K:x_{J,K}\in\mathscr{L}_0(J) \mbox{ is irreducible}\} . 
\end{align} 
Consequently, using \eqref{card_irreducible_antecedents}:
\begin{equation}
\label{eq23s}
a_j  \leq \frac{(2^Q-1)^{1/p}}{j^\beta} \Big(1+ (2^Q-1)\sum_{J=1}^j  2^{ Q(J-1) } 2^{-QJ} \Big)^{1/p}  \leq  \frac{(2^Q-1)^{2/p}}{j^\beta} (1+ j 2^{-Q}) ^{1/p}. 
\end{equation}
Finally, the choice of $\beta$ in our statement provides $a_j \leq  {(2^Q-1)^{2/p}}/{ j^{{2}/{q}}}$ ;
thus the sequence $(a_j)_{j\geq 1}$ belongs to $\ell^q$ and $F\in B^s_{p,q}(\Hn)$.

\mk

In order to compute the multifractal spectrum of $F$ one uses the following lemma.

\begin{lemma}
\label{lem:calculexposant}
For each $x\in [0,1]^3$,  $\displaystyle h_F(x) = s -  \frac{Q}{p}+ \frac  {Q}{p\xi_x}.$  
\end{lemma}
Here   $\xi_x=\xi_x(\mathbb{N})$ is the approximation rate of $x$ by all the dyadic elements.
Assume for a while that Lemma  \ref{lem:calculexposant} holds true ; let us explain how to  conclude from there. Since $\xi_x\in [1,+\infty]$ for each $x$ one first observes that $h_F(x) $  belongs necessarily to the interval $[s- Q/p,s]$. If $h\in (s- Q/p,s]$,  one can observe further that
$$ E_F(h) = \{x\in [0,1]^3 \, : \, h_f(x)=h\} = \left\{x\in [0,1]^3\, : \,  \xi_x =  \frac{Q}{ph - ps +Q} \right\}.$$
Applying \eqref{equaldim}, one deduces that
$$d_F(h) = \dim_H E_F(h) =  \dim_H \widetilde S_{\frac{Q}{ph - ps +Q} } = \frac{Q}{\:\:\frac{Q}{ph - ps +Q} \:\:} = {ph - ps +Q}$$
which is the expected result. If $h =s- Q/p$, the   dimension cannot exceed 0 by the general upper bound given by Theorem~\ref{majspectrumBesov}.
The dimension is exactly 0 because one can find $x\in\Hn$ such that $\xi_x=+\infty$. Such points $x$ are   analogues   in  $\Hn$ of Liouville numbers in $\mathbb{R}$, which are the irrational real numbers that are the ``closest'' to the rationals \cite{HW}.
\cqfd

\mk
\Proof[of Lemma  \ref{lem:calculexposant}]
Consider $x\in \zu^3$ with  $1\leq \xi_x <+\infty$. By definition, for every $\ep>0$, one has the following properties:
\begin{itemize}
\item[(i)]
There exists $J_x>0$ such that for every $j\geq J_x$, for every $k$, $$\delta(x,x_{j,k}) \geq 2^{-j(\xi_x+\ep)}.$$
\item[(ii)]
There exists a strictly increasing sequence of integers $(j_n)_{n\geq 1}$ and a sequence $(k_n)_{n\geq 1} \in \mathcal{Z}^{\mathbb{N}}$ such that  $2^{-j_n}\circ k_n$ is irreducible and $$\delta(x,x_{j_n,k_n})  \leq 2^{-j_n(\xi_x-\ep)}.$$ When $\xi_x=1$ one may take $\ep=0$ in the last inequality.
\end{itemize}

\sk
 To get the lower bound for the H\"older exponent,  consider  dyadic elements $x_{j,k}$ such that   their associated irreducible element $x_{J,K}$ satisfy $J \geq J_x$.  
By item (i), one necessarily has 
$\delta(x,x_{j,k}) = \delta(x,x_{J,K}) \geq 2^{-J (\xi_x+\ep)}.$
  By using that  $2^{-j}$ and $\delta(x,x_{j,k})$ are bounded by above by their sum $ 2^{-j}+\delta(x,x_{j,k})$, we get that
 \begin{eqnarray*}
F^\ep_{j,k} &  =  &  \frac{1}{j^\beta}  \ 2^{ -j(s- Q/p )- J Q/p } \leq ( 2^{-j}+\delta(x,x_{j,k}))^{s- Q/p} (\delta(x,x_{j,k})^{ \frac{1}{\xi_x+\ep}  } )^{    Q/p } \\
 & \leq& (2^{-j}+ \delta(x,x_{j,k}))^{s- Q/p+ \frac{Q}{p(\xi_x+\ep)}} . 
\end{eqnarray*}
This is equivalent to \eqref{PointRegWav}, hence $h_F(x) \geq    s- \frac{Q}{p}+ \frac{Q}{p(\xi_x+\ep)}$.  Letting $\ep$ tend to zero yields the lower bound in Lemma  \ref{lem:calculexposant}.

\mk

Let us  bound by above the H\"older  exponent of $F$ at $x$, by using item $(ii)$. Assume that $1<\xi_x<+\infty$ and fix $\ep>0$ such that $\xi_x-\ep>1$.  For any integer $n\geq 1$, let $\widetilde j_n = [j_n(\xi_x-\ep)]$. Consider the unique dyadic element $x_{ \widetilde j_n , \widetilde k_n}$ such that  $x_{ \widetilde j_n , \widetilde k_n}=x_{  j_n ,  k_n}$. Using  that $2^{-j_n}\circ k_n$ is irreducible, one sees that  
$$F^\ep_{\widetilde j_n, \widetilde k_n} = \frac{1}{(\widetilde j_n)^\beta}2^{- \widetilde j_n(s- Q/p)-j_n Q/p}  \geq  \frac{1}{(\xi_x-\ep)^\beta j_n^\beta}2^{-j_n(\xi_x-\ep) \big(s- \frac{Q}{p}+\frac{Q}{p(\xi_x-\ep)}\big) } .$$
 Hence,  since $ \log( j_n(\xi_x-\ep))$ is negligible with respect to ${ j_n}$ when  $n\to +\infty$, one has
\begin{eqnarray*}
F^\ep_{\widetilde j_n, \widetilde k_n}     \geq      2^{-j_n(\xi_x-\ep)   ({s-   \frac{Q}{p} + \frac{Q}{p(\xi_x-\ep)}  +\ep)}} \geq  d(x,x_{j_n,k_n}) ^{   s-   \frac{Q}{p} + \frac{Q}{p(\xi_x-\ep)}  +\ep}=d(x,x_{\widetilde j_n, \widetilde k_n}) ^{  {s-   \frac{Q}{p} + \frac {Q}{p(\xi_x-\ep)}+\ep} }. 
 \end{eqnarray*}
 This proves that $h_ F (x) \leq   s-   \frac{Q}{p}+ \frac{Q}{p(\xi_x-\ep)}+\varepsilon$, for every $\ep>0$. Letting $\ep$ tend to zero yields the upper bound in Lemma  \ref{lem:calculexposant}. The cases $\xi_x=1$ and $\xi_x=+\infty$ are dealt with similarly.
 \cqfd

 \subsection{The residual set in $B^s_{p,q}( \Hn)$ }
 \label{step3}

Let us define the wavelet version of local spaces. Recall that the set of indices  $\mathscr{L}_0(j)$ was defined in \eqref{defL0}.

\begin{dfn}
The space $B^s_{p,q}([0,1)^3)$ is the closed subspace of $B^s_{p,q}(\Hn)$ defined by
\begin{equation*}
k\notin\mathscr{L}_0(j) \quad \Longrightarrow\quad d_{j,k}^\ep(f)=0.
\end{equation*}
It is  equipped with norm  $
\|f\|_{B^s_{p,q}(\Hn)} =\|f\|_{\infty}+ \|(a_j)_{j\geq 1}\|_{l^q},
$
where $a_j$ are the Besov coefficients of $f$ given  by \eqref{defbesov}.
\end{dfn}

Since   $B^s_{p,q}(\zu^3)$ is separable, one can consider   a countable   sequence $(f_n)_{n\geq 1}$ dense  in  $B^s_{p,q}(\zu^3)$. Let us consider the sequence $(g_n)_{n\geq 1}$ built as follows.
 
 \begin{dfn}For every $n\geq 1$, the wavelet coef\-ficients of $g_n$ up to the generation $j=n-1$ are those of $f_n$ ; for $j\geq n$, the wavelet coefficients of generation $j$ of $g_n$ are those of the function $F$, which are prescribed by equation \eqref{deffjk}.     \end{dfn}

 Since $\|f_n-g_n\|_{B^s_{p,q}(\Hn)}$    tends to zero when $n\to+\infty$,  $(g_n)_{n\geq 1}$ is also dense in $B^s_{p,q}(\zu^3)$.  
 \begin{dfn}
  Let $r_n = \displaystyle  n^{-\beta} 2^{-nQ/p} /2$ with $ \beta $ given by \eqref{deffjk}.  One defines the set $ \widetilde{ \mathcal{R}}$
   $$ \widetilde {\mathcal{R} }= \bigcap_{N \geq 1} \ \bigcup_{n \geq N} \mathscr{B}(g_n,r_n)$$
  where $\mathscr{B}(g,r)=\{f\in B^s_{p,q}([0,1)^3) \, :\, \|f-g\|_{B^s_{p,q}(\Hn)}<r\}$.
   \end{dfn}
The set $ \widetilde {\mathcal{R} } $ is   an intersection of dense open set, hence a residual set in  $B^s_{p,q}(\zu^3)$.
 The   choice for the radius $r_n$ is small enough to ensure that any function $f$ in  $\mathscr{B}(g_n,r_n)$ has its wavelet coefficients at generation $n$ close to those of $g_n$ (and thus to those of $F$).
  
\begin{lemma}
\label{lemfin1}
If $f\in \mathcal{B}(g_n,r_n)$, then  $\displaystyle  |d^\varepsilon_{n,k}(f)-d^\varepsilon_{n,k}(g_n)|  \geq  {|d^\varepsilon_{n,k}(g_n)|}/ {2}\cdotp$
\end{lemma}
\Proof
By definition, one has $d^\varepsilon_{n,k}(g_n) = F_{n,k}^\ep$, $\forall \, k$. Hence, by definition of the Besov norm and the inclusion $\ell^q\subset \ell^\infty$:
$$\left( \sum_{k} \    2^{pn\left(s- Q/p\right)}    |d^\varepsilon_{n,k}(f)-F_{n,k}^\ep|  ^{p} \right)^{1/p} < r_n.$$
In particular,    for any $\ep$ and $k$, 
$$|d^\varepsilon_{n,k}(f)-F_{n,k}^\ep|    \leq  r_n  2^{-n \left(s- Q/p\right)}  \leq  =   2^{-ns}  n^{-\beta}/2.$$
The inequality $J\leq j$ in \eqref{deffjk} reads   $|F_{j,k}^\ep| \geq {2^{-js}}/{j^{\beta}}\cdotp$
Combining both inequalities ensures the result.
\cqfd

\begin{lemma}
If $f\in \widetilde{\mathcal{R}}$,  then  its multifractal spectrum $d_f$ satisfies \eqref{BesovSpectrumequality}.
\end{lemma}
\Proof
Given a function $f\in \widetilde{\mathcal{R}}$, there exists a strictly increasing sequence $ (n_m)_{m\geq 1} $ of integers such that $f \in \mathscr{B}(g_{   n_m},r_{   n_m})$.
Lemma \ref{lemfin1} provides a precise estimate of the wavelet coefficients of $f$, namely for any $m\geq1$:
$$ \frac{1}{2} F_{ {n_m} ,k}^\ep \leq  |d^\varepsilon_{ {n_m} ,k}(f)|  \leq  \frac{3}{2} F_{ {n_m} ,k}^\ep.$$
The same proof as the one developed for Lemma \ref{lem:calculexposant} ensures that 
for any $x\in [0,1)^3$:
$$s- Q/p \leq h_f(x) \leq s -  Q/p +  {Q}/{(p\xi_x(\mathcal{J})} ) \leq s,$$
where $\xi_x(\mathcal{J})$ is the approximation rate by the family $\mathcal{J} = ({n_m})_{m\geq 1} $.
 
\sk
Given $h\in [s- Q/p,s]$ and the unique  $\xi$ such that $h=s- Q/p+ {Q}/{(p\xi)}$, one introduces the set (see Definition \ref{deftaux} and Lemma \ref{leminf}):
$$ \mathcal{E} = S_\xi  (\mathcal{J}) \setminus \bigcup_{i=1}^{+\infty} E^{\leq} _f\left(h- {1}/{i}  \right). $$
By \eqref{maj12} one knows that  $\dim_H  E_f^{\leq}(h') \leq p(h'-s- Q/p)$ for any $h'<h$.
In particular, for every $i\geq 1$,  one has:
$$\dim_H E^{\leq} _f\left(h- {1}/{i}  \right) \leq  p\left(h - {1}/{i} -s- Q/p \right) < p\left(h-s- Q/p \right) =  {Q}/{\xi}\cdotp$$
But according to \eqref{inegmeasure}, one has $\mathcal{H} ^{Q/\xi}(S_\xi( \mathcal{J}) )=+\infty$,
thus $\mathcal{H} ^{Q/\xi}(\mathcal{E}) =+\infty$ and $$\dim_H \, \mathcal{E} \geq  {Q}/ {\xi}\cdotp$$
Next, one observes that  $\mathcal{E} \subset E_f(h)$, since every $x\in  S_\xi  (\mathcal{J})$ satisfies $h_f(x)\leq  h$ and, by definition, $\mathcal{E}$ does not contains those elements $x$ which have
a local exponent strictly smaller than $h$. One can thus finally infer that:
$$\dim_H \, E_f(h) \geq \dim_H \, \mathcal{E} \geq {Q}/ {\xi} = p\left(h-s- Q/p\right) .$$
The converse inequality is provided by Theorem~\ref{majspectrumBesov} because $f\in B^s_{p,q}(\zu^3)$. 
Consequently, the identity \eqref{BesovSpectrumequality} is satisfied.
\cqfd

\mk
To conclude the proof of Theorem \ref{th:besovps}, let us go back to the initial remarks of \S\ref{section:generic_besov}.
The subset $\mathcal{R}_0$ of $B^s_{p,q}(\Hn)$ whose (wavelet) restriction to $[0,1)^3$  satisfies~\eqref{BesovSpectrumequality} and is generic in $B^s_{p,q}([0,1)^3)$ is simply
$$\mathcal{R}_0 = \pi^{-1}(\widetilde{\mathcal{R}})$$
where $\pi : B^s_{p,q}(\Hn)\to B^s_{p,q}(\zu^3)$ is the projection defined in wavelet coefficients by $d_{j,k}^\varepsilon (\pi(f)) = 
d_{j,k}(f)  \cdot {\bf 1\!\!\!1}_{\mathscr{L}_0(j)}(k) $, ${\bf 1\!\!\!1}_{A}(x)$ being equal to 1 if $x\in A$, 0 otherwise.

%
%
%


\section{Generalization to stratified nilpotent groups}
\label{secappendix}

A Carnot group $G$ is a connected, simply connected and nilpotent Lie group whose Lie algebra $\mathfrak{g}$ admits a stratification, \textit{i.e.} for some integer $ {N_G} \geq 1$,
$$\mathfrak{g} = \bigoplus_{k=1}^{ {N_G}} \mathfrak{n}_{k}\qquad\text{where}\qquad [\mathfrak{n}_1,\mathfrak{n}_k]=\mathfrak{n}_{k+1}$$
with $\mathfrak{n}_{ {N_G}}\neq\{0\}$ but $\mathfrak{n}_{ {N_G}+1}=\{0\}$. Let us denote the dimensions $q_k=\operatorname{dim}\mathfrak{n}_k$,
\begin{equation}\label{HomDim}
d=\sum_{k=1}^ {N_G} q_k \qquad\text{and}\qquad Q_G=\sum_{k=1}^ {N_G} kq_k.
\end{equation}
Given a basis $(X_i)_{i=1,\ldots,d}$ of $\mathfrak{g}$ adapted to the stratification, each index $i\in\{1,\ldots,d\}$ can be associated to a unique
$\sigma_i=j\in\{1,\ldots, {N_G}\}$ such that $X_i\in\mathfrak{n}_j$.

Similarly to  \eqref{nabla_Hn}, the horizontal derivatives are the derivatives of the first layer:
$$\nabla_G f=(X_1f,\ldots,X_{q_1}f).$$
The stratification hypothesis ensures that each derivative $X_if$ can be expressed as at most $\sigma_i-1$ commutators of horizontal derivatives.

A Carnot group is naturally endowed with a family of algebra homomorphisms called dilations $\{D_\lambda\}_{\lambda>0}$ that are defined by:
$$\forall i\in\{1,\ldots,d\},\qquad D_\lambda(X_i) = \lambda^{\sigma_i}X_i.$$

The exponential map $\exp:\mathfrak{g}\to G$ is a global analytic diffeomorphism and one can identify $G$ and $\mathfrak{g}$ equipped with the (non commutative if $N\geq2$) law:
$$X\ast Y = \exp^{-1}(\exp(X)\cdot \exp(Y)).$$
 Finally, {as it was the case for $\Hn$}, one can identify $\mathfrak{g}$ to $\R^d$ through the basis $(X_i)_{i=1,\ldots,d}$.
The gauge distance is then defined by
$$\delta(x,y)=\norme[G]{x^{-1}\ast y}, \ \ \  \mbox{ where } \ \ \norme[G]{x}= \Big(\sum_{i=1}^d |x_i|^{2\sigma/\sigma_i}\Big)^{1/2\sigma}$$
and $\sigma=\operatorname{lcm}\{\sigma_1,\ldots,\sigma_d\}$.
This distance is left-invariant and homogeneous of degree 1 with respect to the dilations.  {The triangular inequality stated in Proposition~\ref{TRIANGULAR}
and Corollary~\ref{CorTRIANGULAR} still hold.}

\mk
{With those identifications}, the Haar measure $\ell_G$ on $G$ coincides with the Lebesgue measure on $\R^d$ and
the volume of the ball $\mathcal{B}(x,r)$ is $r^{Q_G}\operatorname{Vol}\mathcal{B}(0,1)$, where $Q_G\geq d$ is defined by~\eqref{HomDim}.
Hausdorff  measures can be defined in a similar fashion to $\Hn$, and   
the  Hausdorff dimension of $G$ is $Q_G $. We refer   to \cite{RigotHausdorff} for further references.

\mk

One can wonder whether the results of multifractal analysis obtained in the present paper still hold in any Carnot group.
We claim that the answer is positive.
Let us list the modifications that are necessary to deal with a Carnot group $G$.

\subsection{Families of wavelets}

Wavelets have been constructed by Lemari\'e \cite{lemarie89} on any Carnot group, but these wavelets have the inconvenient that the number of mother wavelets to be used $\Psi^\ep$ may be infinite, and that the discrete lattice $\mathcal{Z}=\Z^d$ may not be a subgroup any more (which complicates the notion of decomposition of a function on the wavelet basis, and the  analysis of wavelet coefficients as well).

One is naturally led to the wavelet construction proposed by F\"uhr and Mayeli in \cite{FuhrMayeli}. Let us recall the part of their results adapted to our context.

\begin{dfn}
Let $G$ be a Carnot group. A discrete subset $\Gamma\subset G$ is a regular sampling set if there exists a relatively compact Borel set $W\subset G$, neighborhood of the identity, satisfying
\begin{equation}
\label{cov}
G = \bigcup_{\gamma\in \Gamma} \gamma W
\end{equation}
where the equality holds up to a set of $\ell_G$-measure $0$, and such that there is almost no covering in this union, i.e. for all $\alpha\neq \gamma\in \Gamma$, $\ell_G(\alpha W\cap\gamma W)=0$.
\end{dfn}

The role of the lattice $\mathcal{Z}$ in $\Hn$ is now played by $\Gamma$ in $G$.

\begin{dfn}
 For every function $\Psi$, every $ j\in \mathbb{Z}$ and every $\gamma\in \Gamma$, we set
 $$\Psi_{j,\gamma} (x) = \Psi (\gamma^{-1} \ast D_{2^{j}}x) = \Psi(D_{2^j} (x_{j,\gamma}^{-1}\ast x))$$
 {with $x_{j,\gamma}=D_{2^{-j}}\gamma$}.
For every function $f\in B ^s_{p,q}(G) $, the wavelet coefficients of $f$ are
\begin{equation}\label{WAVELETS3}
 d _{j,\gamma}(f) = 2^{jQ_G}\int_{G} f(x)\Psi_{j,\gamma}(x)dx.
\end{equation}
\end{dfn}

Existence of admissible wavelets (i.e.~such that any {smooth enough} function can be reconstructed from its wavelet coefficients) belonging to the Schwartz class on $G$ and having vanishing moments of arbitrary order  is proved for instance in Theorem 4.2 in~\cite{FuhrMayeli}.

\sk
In particular, when adaptating our proofs to general Carnot groups, {one chooses $\psi$ such that} the estimates of the tail of the wavelets \eqref{DECAYPsi} and the vanishing moments \eqref{vanishingmoments} remain unchanged.

\subsection{Taylor polynomials}
Taylor polynomials and estimates of the error term in a Taylor expansion (Theorem \ref{taylor}) hold  on stratified groups (see \cite{follandstein}). 
For general homogeneous groups, a weaker estimate is given in~\cite{follandstein} and explicit Taylor formulas with various remainder terms can also be found in~\cite{TaylorBonfiglioli}.
{In particular, formula \eqref{TAYLOR} above remains valid to compute the  (right) Taylor polynomial of order $N\geq 1 $ at $x_0\in G$ on any Carnot group~$G$~:
$$P_{x_0}(y)  = \sum_{|\alpha|=  |\beta|\leq N }  c_{\alpha,\beta} \nabla_G^\beta f(x_0) y^\alpha.$$
The formula only involves ``horizontal'' derivatives through $\nabla_G f=(X_1f,\ldots,X_{q_1}f)$ but contains homogeneous monomials $y_1^{\alpha_1}\ldots y_d^{\alpha_d}$ of all the coordinates. Homogenity of a multi-index is defined by the proper weights $|\alpha|=\sum_{i=1}^d  \sigma_i\alpha_i$}.

\subsection{Hausdorff dimension}   Throughout the generalization, one   substitutes the new value of the homogeneous dimension $Q_G$. The numerical value of the constants related to the numeration
of neighboring cubes or balls will also have to be modified. Apart from this, definitions and methods are the same.

\subsection{Relations with Besov spaces}

A characterization  of Besov spaces in $G$  in terms of (discrete) wavelet coefficients is similar to the one we used, except that the regular sample is not $\mathcal{Z}=\mathbb{Z}^3$ any more, but rather the regular sampling $\Gamma$.  This is a consequence of Theorems 5.4, 6.1 and 6.7  of \cite{FuhrMayeli}, which can be restated in the following form which suits our context.

\begin{thm}[{\cite{FuhrMayeli}}]
There exists an admissible  wavelet $\Psi$ belonging to the Schwartz class of~$G$ and having infinitely many vanishing moments, and  a regular sampling set $\Gamma$, such that
\begin{equation}
\label{eq25s}
f\in B ^s_{p,q}(G)   \  \ \Rightarrow   \ \ \sum_{j\in \mathbb{Z}}  \left( 2^{j (ps-Q_G)} \left( \sum_{\gamma\in\Gamma} |d _{j,\gamma}(f)|^p \right)\right)^{q/p} <+\infty.
\end{equation}
Reciprocally, if a sequence of coefficients $(c_{j,\gamma})_{j\in \mathbb{Z}, \gamma\in \Gamma}$ satisfies
\begin{equation}
\label{eq12s}
\sum_{j\in \mathbb{Z}}  \left( 2^{j (ps-Q_G)} \left( \sum_{\gamma\in\Gamma} |c _{j,\gamma} |^p \right)\right)^{q/p} <+\infty,
\end{equation}
then the function
$$
f = \sum_{j\in \mathbb{Z}}   \sum_{\gamma\in\Gamma} c _{j,\gamma} \Psi_{j,\gamma}
$$
belongs to $ B ^s_{p,q}(G)$, and the norm $\|f\|_{B^s_{p,q}(G)}$ is equivalent to the sum \eqref{eq12s}.

\sk

Finally, there exists another admissible wavelet  $\widetilde \Psi$ in the Schwartz class of $G$, called {\em dual to $\Psi$}, (depending on the Besov space $B^s_{p,q}(G)$) such that any function $f\in B^s_{p,q}(G)$ can be decomposed as
\begin{equation}
\label{eq24s}
f = \sum_{j\in \mathbb{Z}}   \sum_{\gamma\in\Gamma} \widetilde d_{j,\gamma}(f)  \Psi_{j,\gamma}
 \ \ \ \ \mbox{ with } \ 
 \widetilde d _{j,\gamma}(f) = 2^{jQ}\int_{\Hn} f(x)\widetilde \Psi_{j,\gamma}(x)dx.
 \end{equation}

\end{thm}

The presence of the pair of bi-orthogonal wavelets $(\Psi,\widetilde \Psi)$ implies that the coefficients involved in equation \eqref{eq25s} can   either  be $d_{j,\gamma}(f)$  or $\widetilde d_{j,\gamma}(f)$. Since $\Psi$ and $\widetilde \Psi$ enjoy exactly the same regularity properties, we replace the notation $\widetilde d_{j,\gamma}(f)$ in \eqref{eq24s} by   $d_{j,\gamma}(f)$, by a slight abuse of notations.

In particular, all the methods we developed can easily be adapted using the wavelet coefficients $  d_{j,\gamma}(f)$ for all $j\geq 1$ and $\gamma\in \Gamma$ instead of the family $d_{j,k}(f)$, $j\geq 1$ and $k\in \mathcal{Z}$. {Let us now quickly review the adaptations of each proof.}

\subsection{Results about H\"older regularity}

{
\begin{thm}
Theorems  \ref{thglobal},  \ref{thlocal} and \ref{thgenericCalpha} remain valid on any Carnot group.
\end{thm}
}

\sk
\paragraph{\it Proof of Theorem \ref{thglobal}.}
In section \ref{sec31}, the only modification consists in naming $\vartheta$ a solution of $\mathcal{L}^M\vartheta = \psi$ for some arbitrarily large integer $M$ and where
$\mathcal{L}=-(X_1^2+\ldots+X_{q_1}^2)$ denotes the hypoelliptic Laplace operator on $G$. Thanks to \cite{FuhrMayeli}, one can choose $\psi$ properly to ensure that $\vartheta$ has at least one vanishing moment and fast decay at infinity, which is all that is required for the proof to work. When $[s]$ is odd, the last integration by part with respect to $(X,Y)$ is obviously replaced by
and integration by part against each $(X_1,\ldots,X_{q_1})$ and produces $q_1$ terms, all dealt with in a similar way.
Section \ref{sec32} remains exactly unchanged.

\sk
\paragraph{\it Proof of Theorem \ref{thlocal}.}
Section \ref{sec41} also remains unchanged because, as noticed above, Corollary~\ref{CorTRIANGULAR} is still valid on $G$. The conclusive statement~\eqref{PointRegWav} should obviously read:
$$\delta(x_{j,\gamma},x_0)<R \quad\Longrightarrow\quad |d_{j,\gamma}(f)|\leq C 2^{-js}\left(1+2^j \delta(x_{j,\gamma},x_0)\right)^s$$
In section \ref{sec42}, the reconstruction  of $f$ from its wavelet coefficients to be used is \eqref{eq24s}, and rewrites
$$f(x)=f^\flat(x)+\sum_{j=1}^\infty \sum_{\gamma\in\Gamma} d_{j,\gamma}(f) \Psi_{j,\gamma}(x),$$
with $f^\flat(x)$ being a smooth function.

The polynomial $P$ suitable for the pointwise H\"older estimate is $P=P^\flat + \sum_{j=1}^\infty P_j$
where $P_j$ is the Taylor expansion of $\sum_{\gamma\in\Gamma} d_{j,\gamma}(f) \Psi_{j,\gamma}$. As noticed above,
the Taylor expansion formula~\eqref{TAYLOR} remains valid on $G$ so the rest of the section is unchanged.

\sk
\paragraph{\it Proof of Theorem \ref{thgenericCalpha}.}
Section \ref{genericcalpha} is an abstract game of seeking wavelets coefficients of the proper order and rounding them up to the closest dyadic integer. It only connects to the ambient space through the application of Theorems  \ref{thglobal} and \ref{thlocal} that we now know to hold true on $G$. The sole modification is the notation $d_{j,\gamma}(f)$ instead of $d^\varepsilon_{j,k}(f)$.
\cqfd

\subsection{Results about Besov spaces and diophantine approximation in $G$.}

One powerful property of the mass transference principle by Beresnevich, Dickinson and Velani \cite{BDV} and similar results in heterogeneous situations \cite{BS,BS2,BS3} is that these theorems
not only apply to approximation by dyadics or rationals in Euclidian settings but also to all sufficiently well-distributed systems of points in doubling metric spaces.

\sk
 The definition of the regular sampling $\Gamma$ and its associated tile $W$ such that \eqref{cov} holds true implies the two following properties:
 \begin{enumerate}
 \item[(C1)]
 Since $W$ is compact, there exists a sufficiently large $M_G>0$ such that any ball $\mathcal{B}(x,M_G)$ contains at least one point $\gamma\in \Gamma$.
 \item[(C2)]
 Since the union \eqref{cov} is constituted by sets whose intersections are always of $\ell_G$-measure 0 and $W$ is bounded,  there exists another constant $N_G>0$ such that for  every $x\in G$, the ball $\mathcal{B}(x,M_G)$ contains at most $N_G$ points belonging to  $\Gamma$.
 \end{enumerate}
 These properties are analogues in $G$  to Lemma \ref{number cover} in the Heisenberg group $\Hn$.  
 One concludes that 
$$G = \bigcup_{\gamma\in \Gamma} \mathcal{B}(\gamma,M_G),$$
and that there is almost no   redundancy in the covering, \textit{i.e.}~for every $x\in G$, the cardinality of those $\gamma\in \Gamma$ such that $x\in  \mathcal{B}(\gamma,M_G)$ is bounded from
above by $N_G$ uniformly in $x\in G$. Immediately, one also deduces an analogue of~\eqref{cover}:
$${W} =  \limsup_{j\to +\infty} \bigcup_{\gamma\in \mathcal{L}_0(j_m)} \mathcal{B}(x_{j,\gamma},2^{-j}M_G)$$
{where $\mathcal{L}_0(j)=\{\gamma\in\Gamma: x_{j,\gamma}=D_{2^{-j}}\gamma \in W\}$. The compact tile $W$ is a natural candidate to replace $[0,1)^3$ on $\Hn$}.
The notion of approximation rate and the sets $\mathcal{S}_\xi(\mathcal{J})$ and $\widetilde {\mathcal{S}}_\xi(\mathcal{J})$ are perfectly defined (recall Definition \ref{deftaux}), and have the same interpretation as in the case of~$\Hn$.

\sk
{
We are now ready to state our last result.
\begin{thm} Theorems \ref{majspectrumBesov} and \ref{th:besovps} remain valid on any Carnot group.
\end{thm}
}

\sk
\paragraph{\it Proof of Theorem \ref{majspectrumBesov}.}

A careful reading of Section \ref{secbesov1} shows that the arguments go through by simply replacing $\mathcal{Z}$ by $\Gamma$.
Indeed, Lemma \ref{lemma1} is a general counting argument for convergent series
and Lemma \ref{leminf} requires only counting  and coverings arguments that are exactly items (C1) and (C2)  explained a few lines above. 
One deduces that every functions $f \in B^s_{p,q}(G)$ satisfies
$$\forall \ h\geq s-Q_G/p, \ \ \ d_f(h) \leq \min(Q_G, ph-ps+Q_G),$$
as  was the case for $\Hn$.
\sk 

\paragraph{\it Proof of Theorem \ref{th:besovps}.}

{
Let us carefully go through Section~\ref{section:generic_besov}. As mentioned above, $W$ replaces $[0,1)^3$.
The first change concerns  the cardinality of $\mathcal{L}_0(j)$.
As $W$ is a neighborhood of the origin, there exists $\varepsilon>0$ such that $\mathcal{B}(\operatorname{Id},\varepsilon)\subset W$. By~\eqref{TRIANGULAR}, there also exists
a constant $C$ such that $\ell_G(\bigcup_{w\in W} \mathcal{B}(w,\varepsilon 2^{-j})) \leq C$ uniformly in $j\geq 0$. If $\gamma_1\neq \gamma_2\in \mathcal{L}_0(j)$ then
$\ell_G \big ( \mathcal{B}(x_{j,\gamma_1},\varepsilon 2^{-j}) \cap \mathcal{B}(x_{j,\gamma_2},\varepsilon 2^{-j})\big ) = 2^{-jQ_G} \ell_G \big ( \mathcal{B}(\gamma_1,\varepsilon) \cap \mathcal{B}(\gamma_2,\varepsilon) \big) $ which is zero
because $\mathcal{B}(\gamma_i,\varepsilon)\subset \gamma_i W$ and $\ell_G(\gamma_1 W\cap \gamma_2 W)=0$. Therefore
\begin{equation}\label{card_L0jBis}
\#\mathcal{L}_0(j) \leq \frac{C}{\ell_G\big(\mathcal{B}((\operatorname{Id},\varepsilon 2^{-j}) \big)} \leq \widetilde{C} 2^{j Q_G}
\end{equation}
which replaces~\eqref{card_L0j}.
}

\sk
In order to prove the optimality of the upper bound for the multifractal spectrum of functions in $B^s_{p,q}(G)$, an ``optimal'' function $F$ was built in Proposition \ref{propF}. Here, the new function to be studied is called $F_G$ and is defined as the sum 
$$F_G = \sum_{j\in \mathbb{Z}}   \sum_{\gamma\in\Gamma} F ^G_{j,\gamma} \Psi_{j,\gamma}$$
where the wavelet coefficients are
\begin{equation}
\label{deffjk2}
F^G_{j,\gamma} :=   \begin{cases}\dfrac{2^{ -j(s- Q_G/p )- J Q_G/p }}{j^\beta} & \mbox{if }{x_{j,\gamma} \in W}\\
0 & \mbox{otherwise,}
\end{cases}
\end{equation}
where $J$ is the minimal positive integer such that $D_{2^J}x_{j,\gamma}=D_{2^{J-j}}\gamma \in \Gamma$. This naturally replaces the notion of irreducibility of dyadics given in Definition \ref{defirred}:
{
for any positive integer~$j'$ such that $\gamma'=D_{2^{j'}} \gamma\in\Gamma$, one has $x_{j,\gamma}=x_{j+j',\gamma'}$, thus $x_{j,\gamma}$ is irreducible if it cannot be written as $x_{j'',\gamma''}$ with $0\leq j''<j$.
One checks easily that $F_G$ belongs to $B^s_{p,q}(G)$. The main estimate that replaces \eqref{eq22s} and \eqref{eq23s} are
$$\|2^{j(s-Q_G/p)}F_{j,\gamma}\|_{\ell^p(\Gamma)} \leq \frac{1}{j^{2/q}}\left(\frac{1}{j}\sum_{J=0}^j 2^{-Q_GJ}
\#\{\gamma\in \mathcal{L}_0(J) \,;\, x_{J,\gamma} \text{ irreducible}\} \right)^{1/p} \!\leq \frac{\widetilde{C}^{1/p}}{j^{2/q}}$$
which, as expected, belongs to  $\ell^q(j\in\mathbb{N})$.
The last inequality results from \eqref{card_L0jBis} which is slightly rougher than the right-hand side of~\eqref{card_irreducible_antecedents} but still sufficient for our purpose.
}

\mk

{In the proof of Proposition~\ref{propequaldim}, the constant $\widetilde{C}$ of~\eqref{card_L0jBis} also appears in the upper bound for the Haussdorff pre-measure that now reads
$$\mathcal{H}^d_\eta(S_\xi(\mathcal{J}))\leq C \sum_{j\geq n} (\widetilde{C}2^{jQ_G})(2^{-j\xi})^d \leq C\widetilde{C}  2^{n (Q_G -d \xi)}$$
and still tends to zero as $n\to + \infty$ when $d>Q_G/\xi$.}
Conversely, let us observe that the techniques we used to find lower bounds for the Hausdorff dimensions of sets extend to Carnot groups thus ensuring the second half of Proposition~\ref{propequaldim}.
It is an easy matter to check that the Haar measure $\ell_G$ satisfies the three conditions (H1), (H2) and (H3) where the set   $\{x_{j,k} \,  :\, j\in\mathcal{J}, \enspace k\in\mathcal{Z}\}$ is replaced by the discrete set $\{x_{j,\gamma}=D_{2^{-j} } \gamma\,  :\, j\in\mathcal{J}, \enspace \gamma\in\Gamma\}$. This implies that the mass transference principle (Theorem \ref{thm_BDV}) holds true on $G$ as it did
in the Heisenberg group $\Hn$. Hence, all the arguments developed to find lower bounds for the Hausdorff multifractal spectrum of typical functions in  $B^s_{p,q}(\Hn)$ can be extended
without alteration to the Carnot group $G$ and its Besov space $B^s_{p,q}(G)$.

\sk
{The last alteration consists in defining the wavelet-local space $B^s_{p,q}(W)$ by the criterion
$$x_{j,\gamma}\notin W \quad\Longrightarrow\quad d_{j,\gamma}(f)=0.$$
The rest of the Section~\ref{step3} remains unchanged.
}
\cqfd

\mk
Taking those remarks in consideration, one can assert that Theorems~\ref{thglobal} to~\ref{th:besovps} remain valid on any Carnot group.

\subsection{{Open problems}}
Further generalization (\textit{e.g.} to the realm of homogeneous groups) are not   straightforward.
For example, even though the metric structure of homogeneous groups is still defined in a similar fashion to the gauge distance  on Carnot groups, the notion of horizontal derivative ceases to exist, which changes deeply the nature of the Taylor formula and its remainder  \cite{TaylorBonfiglioli} and thus the subsequent analysis.
The construction and analysis of wavelets in such a general setting is also an active area of mathematics.
 
\sk
The reader might also be interested in the following works concerning wavelets
on compact Lie groups \cite{stein_LP},  on general Lie groups \cite{Triebel}, on homogeneous spaces \cite{EbertWirth} and even riemannian manifolds \cite{GellerMayeli}.



\begin{thebibliography}{2010}

\bibitem{Mat1}
Balogh, Z.;  Durand  Cartagena, E.;  F\"assler, K.;  Mattila P., Tyson,  J., The effect of projections on dimension in the Heisenberg group, to appear in Revista Math. Iberoamericana.

\bibitem{Mat2}
Balogh, Z.;   F\"assler, K.;  Mattila P., Projection and slicing theorems in Heisenberg groups, Adv. Math. 231 (2012),   569-604.


\bibitem{BS}
Barral, J.; Seuret, S.,  Heterogeneous ubiquitous systems in $\Bbb R\sp d$ and Hausdorff dimension. Bull. Braz. Math. Soc. (N.S.) 38 (2007), no. 3, 467-515.


\bibitem{BS2}
Barral, J.; Seuret, S., Ubiquity and large intersections properties under digit frequencies constraints. Math. Proc. Cambridge Philos. Soc. 145 (2008), no. 3, 527-548.


\bibitem{BS3}
Barral, J.;  Seuret, S., A localized Jarnik-Besicovitch theorem. Adv. Math. 226 (2011), no. 4, 3191-3215.

\bibitem{BDV}
Beresnevich, V.; Dickinson, D.; Velani, S.,
Measure theoretic laws for lim sup sets.
Mem. Amer. Math. Soc. 179 (2006), no. 846

\bibitem{BP}
Berhanu, S.; Pesenson, I.,
The trace problem for vector fields satisfying Hörmander's condition.
Math. Z. 231 (1999), no. 1, 103-122.

\bibitem{TaylorBonfiglioli}
Bonfiglioli, A.,
Taylor formula for homogeneous groups and applications.
Math. Z. 262 (2009), no. 2, 255-279. 

\bibitem{BC}
Bony, J.-M.; Chemin, J.-Y.,
Espaces fonctionnels associés au calcul de Weyl-Hörmander.
Bull. Soc. Math. France 122 (1994), no. 1, 77-118.

\bibitem{CCX}
Cancelier, C. E.; Chemin, J.-Y.; Xu, C. J.,
Calcul de Weyl et opérateurs sous-elliptiques.
Ann. Inst. Fourier 43 (1993), no. 4, 1157-1178. 

\bibitem{CX}
Chemin, J.-Y.; Xu, C. J.,
Sobolev embeddings in Weyl-Hörmander calculus. Geometrical optics and related topics (Cortona, 1996), 79-93, 
Progr. Nonlinear Differential Equations Appl., 32, Birkhäuser Boston, Boston, MA, 1997. 
%


\bibitem{Cygan}
Cygan, J., 
Subadditivity of homogeneous norms on certain nilpotent groups.
Proc. Amer. Math. Soc. 83 (1981), 69--70.

\bibitem{DGN}
Danielli, D.; Garofalo, N.; Nhieu, D.-M.,
Non-doubling Ahlfors measures, perimeter measures and the characterization of the trace spaces of Sobolev functions in Carnot-Carathéodory spaces.
Mem. Amer. Math. Soc. 182 (2006), no. 857.

\bibitem{EbertWirth}
Ebert, S.; Wirth, J.,
Diffusive wavelets on groups and homogeneous spaces.
Proc. Roy. Soc. Edinburgh Sect. A 141 (2011), no. 3, 497-520.

\bibitem{folland}
Folland, G. B.
Subelliptic estimates and function spaces on nilpotent Lie groups. 
Ark. Mat. 13 (1975), no. 2, 161-207. 

\bibitem{follandstein}
Folland, G. B.; Stein, E. M.,
Hardy spaces on homogeneous groups. 
Mathematical Notes, 28. Princeton University Press, Princeton, N.J., 1982.

\bibitem{FuhrMayeli}
F\"{u}hr, H.; Mayeli, A.,
Homogeneous Besov spaces on stratified Lie groups and their wavelet characterization.
J. Funct. Spaces Appl. 2012, Art. ID 523586.


\bibitem{Ajoutreferee}
 Furioli, G., 
Melzi, C., Veneruso A., 
 Littlewood-Paley decompositions and Besov spaces on Lie groups of polynomial growth, Math. Nach.  279 (9-10) (2006) 1028--1040.
\bibitem{GellerMayeli}
Geller, D.; Mayeli, A.,
Besov spaces and frames on compact manifolds.
Indiana Univ. Math. J. 58 (2009), no. 5, 2003-2042. 


\bibitem{HW}
Hardy, G. H.; Wright, E.M.,
An Introduction to the Theory of Numbers.
Oxford Univ. Press. 2008.



\bibitem{Jaffard_generic}
Jaffard, S.,
On the Frisch-Parisi conjecture.
J. Math. Pures Appl. (9) 79 (2000), no. 6, 525-552.

\bibitem{jaffleaders}
Jaffard, S.; Lashermes, B.; Abry, P.,
Wavelet leaders in multifractal analysis.
Wavelet analysis and applications, 201-246, 
Appl. Numer. Harmon. Anal., Birkhäuser, Basel, 2007.

\bibitem{Besicovitch2}
Korányi, A.; Reimann, H. M.,
Foundations for the theory of quasiconformal mappings on the Heisenberg group. 
Adv. Math. 111 (1995), no. 1, 1-87. 

\bibitem{lemarie89}
Lemarié, P. G.,
Base d'ondelettes sur les groupes de Lie stratifi\'{e}s.
Bull. Soc. Math. France 117 (1989), no. 2, 211-232. 

\bibitem{MV}
Mustapha, S.; Vigneron, F.,
Construction of Sobolev spaces of fractional order with sub-Riemannian vector fields.
Ann. Inst. Fourier (Grenoble) 57 (2007), no. 4, 1023-1049. 

\bibitem{Besicovitch3}
Rigot, S.,
Counter example to the Besicovitch covering property for some Carnot groups equipped with their Carnot-Carathéodory metric. (English summary) 
Math. Z. 248 (2004), no. 4, 827-848. 

\bibitem{RigotHausdorff}
Rigot, S.,
Isodiametric inequality in Carnot groups.
Ann. Acad. Sci. Fenn. Math. 36 (2011), no. 1, 245-260. 

\bibitem{stein}
Stein, E. M.,
Harmonic analysis: real-variable methods, orthogonality and oscillatory integrals.
Princeton Mathematical Series, 43. Monographs in Harmonic Analysis, III. Princeton University Press, Princeton, NJ, 1993.

\bibitem{stein_LP}
Stein, E. M.,
Topics in harmonic analysis related to the Littlewood-Paley theory. 
Annals of Mathematics Studies, No. 63 Princeton University Press, Princeton, N.J., 1970.

\bibitem{saka}
Saka, K.,
Besov spaces and Sobolev spaces on a nilpotent Lie group. 
Tôhoku Math. J. (2) 31 (1979), no. 4, 383-437. 

\bibitem{Besicovitch1}
Sawyer, E.; Wheeden, R. L.,
Weighted inequalities for fractional integrals on Euclidean and homogeneous spaces. 
Amer. J. Math. 114 (1992), no. 4, 813-874. 

\bibitem{Triebel}
Triebel, H.,
Function spaces on Lie groups, the Riemannian approach. 
J. London Math. Soc. (2) 35 (1987), no. 2, 327-338. 


 \bibitem{Varadarajan} 
Varadarajan, V. S., 
Lie groups, Lie algebras, and their representations.
Graduate Texts in Mathematics, Springer, 1984.

\bibitem{vigneron}
Vigneron, F.,
The trace problem for Sobolev spaces over the Heisenberg group. (English summary) 
J. Anal. Math. 103 (2007), 279-306. 

\end{thebibliography}
\end{document}